\documentclass[11pt]{article}
\usepackage{amsfonts,amssymb,amsmath,mathrsfs,bbm}
\usepackage{color,latexsym,amsfonts}

\numberwithin{equation}{section}
\newtheorem{theorem}{Theorem}[section]
\newtheorem{lemma}{Lemma}

\newtheorem{corollary}[theorem]{Corollary}

\newtheorem{remark}{Remark}

\def\<{\langle}
\def\>{\rangle}

\begin{document}

\bigskip\bigskip \noindent{\Large \textbf{Limit theorems for supercritical
MBPRE with linear fractional offspring distributions }}\footnote{%
\noindent Supported by NSFC (NO.11531001, 11626245), High-end Foreign Experts Recruitment Program (GDW20171100029).}

\noindent {%\normalsize\sf
Wenming Hong\footnote{%
School of Mathematical Sciences \& Laboratory of Mathematics and Complex
Systems, Beijing Normal University, Beijing 100875, P.R. China. Email:
wmhong@bnu.edu.cn} ~ Minzhi Liu\footnote{%
School of Mathematical Sciences \& Laboratory of Mathematics and Complex
Systems, Beijing Normal University, Beijing 100875, P.R. China. Email:
liuminzhi@mail.bnu.edu.cn}} ~ Vladimir Vatutin \footnote{%
Beijing Normal University, Beijing 100875, P.R. China and Steklov
Mathematical Institute Gubkin street 117 996 Moscow GSP- I Russia Email:
vatutin@mi.ras.ru}

\begin{center}
\begin{minipage}{12cm}
\begin{center}\textbf{Abstract}\end{center}
\footnotesize
We investigate the limit behavior of supercritical multitype branching processes in random environments
with linear fractional offspring distributions and show that there exists a phase transition in the behavior
of local probabilites  of the process affected by strongly and intermediately supercritical regimes.
Some conditional limit theorems can also be obtained from the representation of generating functions.

\bigskip

\mbox{}\textbf{Keywords:}\quad multitype branching process in random environments, supercritical, random walk, conditional limit theorem; \\
\mbox{}\textbf{Mathematics Subject Classification}:  Primary 60J80;
secondary 60G50.

\end{minipage}
\end{center}

\section{ Introduction\label{Introd}}

Branching processes in random environments (BPRE), which are one of the
interesting and important generalizations of the Galton-Watson processes,
have been first introduced by Smith and Wilkinson \cite{smith-wilkinson} for
the case of i.i.d. environmental situation. Athreya and Karlin in \cite%
{athreya-karlin1} and \cite{athreya-karlin2} extended this model to a more
general environmental situation. Since then, Weissner \cite{weissner},
Kaplan \cite{kaplan} and Tanny \cite{tanny} considered the multitype
branching processes in random environments (MBPRE) of Smith-Wilkinson model
and Athreya-Karlin model respectively, and have got the extinction criteria
for them. Recently, some interesting properties such as a kind of phase
transition were noted in the subcritical and supercritical cases for BPRE
model in \cite{birkner-geiger-kersting}, \cite{AGKV2}, \cite{ABKV10}, \cite%
{ABKV11}, and \cite{boinghoff}, \cite{HL2012}, \cite{GLM2016}. To be more
specific, it was shown that the survival probability decreases at different
rates in the strongly, intermediately and weakly subcritical phases. These
finding were complemented in \cite{dyakonova08} and \cite{dyakonova13} by
describing the asymptotics of the survival probability for some classes of
subcritical MBPRE.

All these papers deal with the statements concerning the events occurring
with negligible probability when time goes to infinity (survival of
subcritical processes, conditional limit theorems for subcritical processes
given their survival, small population size for single-type supercritical
processes given their survival). However, up to now the supercritical MBPRE
were not investigated in detail. In the present paper, we study the
distribution of the population size for a class of supercritical MBPRE and
investigate the asymptotic behavior of the probabilities of some rare events
including the probability to have a small number of particles in a
supercritical MBPRE. In particular, we generalize some results from \cite%
{boinghoff} concerning the asymptotic behavior of the probabilities of some
rare events for single-type BPRE to the multitype setting.

\section{Preliminaries\label{Sec_Prelim}}

Dealing with multitype processes it is impossible to avoid complicated
notation and formulas. We begin our description with introducing some
standard agreement for $K$-dimensional vectors and matrices. First we make
no distinction in the notation of row and column vectors. It will be clear
from the context which form (row or column) of the respective vector is
used. In addition,

\begin{itemize}
\item we denote by $\mathbf{e}_{j}$ a $K$-dimensional vector whose $j$-th
component equals $1$ and all others equal zero;

\item all zero and all one vectors will be written as $\mathbf{0}=(0,\ldots
,0)$ and $\mathbf{1}=(1,\ldots ,1);$

\item for $\mathbf{x}=\left( x_{1},\ldots ,x_{K}\right) $ and $\mathbf{y}%
=\left( y_{1},\ldots ,y_{K}\right) $ we set
\begin{equation*}
|\mathbf{x}|=\sum\limits_{i=1}^{K}|x_{i}|,\ \left( \mathbf{x},\mathbf{y}%
\right) =\sum\limits_{i=1}^{K}x_{i}y_{i};
\end{equation*}

\item for a positive $K\times K$ matrix $\mathbf{A}=\left( A(i,j)\right)
_{i,j=1}^{K}$, denote its Perron root by $\rho (\mathbf{A})$ and let $%
\mathbf{v}(\mathbf{A})=\left( v_{1}(\mathbf{A}),\ldots ,v_{K}(\mathbf{A}%
)\right) $ and $\mathbf{u}(\mathbf{A})=\left( u_{1}(\mathbf{A}),\ldots
,u_{K}(\mathbf{A})\right) $ be the strictly positive left and right
eigenvectors of $\mathbf{A}$ corresponding to $\rho (\mathbf{A})$ and
meeting the scaling
\begin{equation}
\left( \mathbf{1},\mathbf{u}(\mathbf{A})\right) =1,\quad \left( \mathbf{v}(%
\mathbf{A}),\mathbf{u}(\mathbf{A})\right) =1;  \label{3}
\end{equation}

\item finally, using the notation $\mathbb{N}_{0}:=\{0,1,2,\ldots \}$ and
taking $\mathbf{s}=\left( s_{1},\ldots ,s_{K}\right) \in \lbrack 0,1]^{K}$
and $\mathbf{z}=\left( z_{1},\ldots ,z_{K}\right) \in \mathbb{N}_{0}^{K}$ we
set $\mathbf{s}^{\mathbf{z}}:=\prod\limits_{i=1}^{K}s_{i}^{z_{i}}.$
\end{itemize}

We now introduce notation related to the $K-$type BPRE's.

Let $\mathfrak{P}(\mathbb{N}_{0}^{K})$ be the space of all probability
measures on the set $\mathbb{N}_{0}^{K}$ of $K$-dimensional vectors with
integer-valued nonnegative components. For $f\mathbf{\in }\mathfrak{P}(%
\mathbb{N}_{0}^{K})$ we denote its weights by $f[\mathbf{z}]$, $\mathbf{z=}%
\left( z_{1},...,z_{K}\right) \in \mathbb{N}_{0}^{K}.$ \ We put%
\begin{equation*}
f\mathbf{(s):=}\sum_{\mathbf{z}\in \mathbb{N}_{0}^{K}}f[\mathbf{z}]\mathbf{s}%
^{\mathbf{z}},\quad \mathbf{s=}\left( s_{1},\ldots ,s_{K}\right) \in \lbrack
0,1]^{K}.
\end{equation*}%
The resulting function on the $K$-dimensional cube $[0,1]^{K}$ is the
generating function of the measure $f$. \ We take the liberty here to denote
the measure and its generating function by one and the same symbol $f$. \ We
need to consider $K-$dimensional vectors
\begin{equation*}
\mathbf{f}=\left( f^{(1)},...,f^{(K)}\right) \in \mathfrak{P}(\mathbb{N}%
_{0}^{K})\times \cdots \times \mathfrak{P}(\mathbb{N}_{0}^{K}):=\mathfrak{P}%
^{K}(\mathbb{N}_{0}^{K})
\end{equation*}%
of probability measures on $\left( \mathbb{N}_{0}^{K}\right) ^{K}:=\mathbb{N}%
_{0}^{K}\times \mathbb{N}_{0}^{K}\times \cdots \times \mathbb{N}_{0}^{K}$
where $f^{(i)}\in \mathfrak{P}(\mathbb{N}_{0}^{K})$ \ for all $i=1,2,...,K$.

Now we specify a $K$-type branching process in varying environment on the
underlying probability space $\left( \Omega ,\mathcal{F},\mathcal{P}\right).
$

\textbf{Definition 1}. A sequence $v=(\mathbf{f}_{1},\mathbf{f}_{2},\ldots )$
of vector-valued probability measures on $\left( \mathbb{N}_{0}^{K}\right)
^{K}$ is called a varying environment.

\textbf{Definition 2}. Let $v=(\mathbf{f}_{n},n\geq 1)$ be a varying
environment. A stochastic process $\mathcal{Z}$ $=$ $\left\{ \mathbf{Z}%
_{n}=(Z_{n}(1),\ldots ,Z_{n}(K)),n\geq 0\right\} $ with values in $\mathbb{N}%
_{0}^{K}$ is called a branching process with varying environment $v$, if for
any $\mathbf{z}\in \mathbb{N}_{0}^{K}$ and $n\geq 1$%
\begin{equation*}
\mathcal{P}\left( \mathbf{Z}_{n}=\mathbf{z\,}|\,\mathbf{Z}_{0},\ldots ,%
\mathbf{Z}_{n-1}\right) =\left( \mathbf{f}_{n}^{\mathbf{Z}_{n-1}}\right) [%
\mathbf{z}].
\end{equation*}

This definition admits the following interpretation in probabilistic terms
for $n\geq 1$: given $\mathbf{Z}_{0},\ldots,\mathbf{Z}_{n-1}$ the random
vector $\mathbf{Z}_{n}=(Z_{n}(1),\ldots ,Z_{n}(K))$ may be realized as the
sum of independent random vectors $\mathbf{Y}_{i,n}^{(j)}=\left(
Y_{i,n}^{(j)}(1),\ldots ,Y_{i,n}^{(j)}(K)\right) $ with distributions $%
f_{n}^{(j)},j=1,\ldots ,K$%
\begin{equation*}
\mathbf{Z}_{n}=\sum_{j=1}^{K}\sum_{i=1}^{Z_{n-1}(j)}\mathbf{Y}_{i,n}^{(j)}.
\end{equation*}%
Thus, informally, $Z_{n}(j)$ is the number of type $j$ individuals of some
population in generation~$n$, where all individuals reproduce independently
of each other and of $\mathbf{Z}_{n}$, and where $f_{n}^{(j)}$ is the
distribution of the offspring vector $\mathbf{Y}_{n}^{(j)}=\left(
Y_{n}^{(j)}(1),\ldots,Y_{n}^{(j)}(K)\right) $ of a type $j$ individual in
generation $n-1$. The distribution of $\mathbf{Z}_{0}$, which is the initial
distribution of the population, may be arbitrary. In the sequel we often
choose $\mathbf{Z}_{0}=\mathbf{e}_{j}$ for some $j\in \left\{
1,\ldots,K\right\} $.

We associate with a vector-valued probability generating function $\mathbf{f}
$ $=\left( f^{(1)},\ldots,f^{(K)}\right) $ the mean matrix
\begin{equation*}
\mathbf{m}=\mathbf{m}\left( \mathbf{f}\right) =(m(i,j))_{i,j=1}^{K}:=\left(
\frac{\partial f^{(i)}}{\partial s_{j}}(\mathbf{1})\right) _{i,j=1}^{K}
\end{equation*}%
and use $\mathbf{m}_{n}$ for $\mathbf{m}\left( \mathbf{f}_{n}\right) $ in
the sequel.

In what follows we write $\mathcal{P}_{\mathbf{z},v}(\mathbf{\cdot })$ and $%
\mathcal{E}_{\mathbf{z},v}\left[ \mathbf{\cdot }\right] $ for the
distribution and expectation of the process in varying environment $v$ with $%
\mathbf{Z}_{0}=\mathbf{z}$. \ If the particular form of the environment is
of no importance we often write $\mathcal{P}_{\mathbf{z}}(\mathbf{\cdot })$
and $\mathcal{E}_{\mathbf{z}}\left[ \mathbf{\cdot }\right] $ for $\mathcal{P}%
_{\mathbf{z},v}(\mathbf{\cdot })$ and $\mathcal{E}_{\mathbf{z},v}\left[
\mathbf{\cdot }\right] .$

Let us pass to the description of $K$-type branching processes in random
environment. To this aim we endow the space $\mathfrak{P}^{K}(\mathbb{N}%
_{0}^{K})$ of vector-valued probability measures on $\left( \mathbb{N}%
_{0}^{K}\right) ^{K}$ with the metric $d_{TV}$ of total variation, given by%
\begin{equation*}
d_{TV}(\mathbf{f,g})=\frac{1}{2K}\sum_{\mathbf{z}\in \mathbb{N}%
_{0}^{K}}\left\vert \mathbf{f}[\mathbf{z}]-\mathbf{g}[\mathbf{z}]\right\vert
,\quad \mathbf{f,g}\in \mathfrak{P}^{K}(\mathbb{N}_{0}^{K})
\end{equation*}%
and with the respective Borel $\sigma $-algebra. This allows us to consider
vector-valued random probability measures on $\left( \mathbb{N}%
_{0}^{K}\right) ^{K}$, which are random vectors $\mathbf{F=}\left(
F^{(1)},...,F^{(K)}\right) $ with values in $\mathfrak{P}^{K}(\mathbb{N}%
_{0}^{K})$ and components specified by probability generating functions
\begin{equation}
F^{(i)}(\mathbf{s}):=\sum_{\mathbf{z}\in \mathbb{N}_{0}^{K}}F^{(i)}[\mathbf{z%
}]\mathbf{s}^{\mathbf{z}},\ i=1,\ldots ,K.  \label{Defff}
\end{equation}

\textbf{Definition 3}. A sequence $\mathcal{V}$ $=\left( \mathbf{F}_{1},%
\mathbf{F}_{2},...\right) $ of random probability measures on $\left(
\mathbb{N}_{0}^{K}\right) ^{K}$ is called a random environment for a $K$%
-type BPRE. An i.i.d. environment is a random environment where $\mathbf{F}%
_{1},\mathbf{F}_{2},\ldots$ are i.i.d. copies of some vector-valued random
probability measure $\mathbf{F}$ on $\left( \mathbb{N}_{0}^{K}\right) ^{K}$.

In the sequel we use the symbols $\mathbb{P}$ and $\mathbb{E}$ for
probabilities and expectations for the branching processes evolving in
random environment in contrast with the symbols $\mathcal{P}$ and $\mathcal{E%
}$ used for the branching processes evolving in varying environment. We
denote by $P$ the law of the environment and by $\mathbb{P}$ the annealed
law of the MBPRE respectively, and use the symbols $E$ and $\mathbb{E}$ for
the corresponding expectations. Then $\mathbb{P}(\cdot ):=\int \mathcal{P}%
_{v}\left( \cdot \right) P(dv)$.

\textbf{Definition 4}. Let $\mathcal{V}$ be a random environment. A process $%
\mathbf{Z}_{n}=(Z_{n}(1),\ldots ,Z_{n}(K))$ with values in $\mathbb{N}%
_{0}^{K}$ is called a $K-$type branching process with random environment $%
\mathcal{V}$ , if for each varying environment $v$ and for all $\mathbf{z},%
\mathbf{z}_{1},\ldots ,\mathbf{z}_{m}\in \mathbb{N}_{0}^{K}$ we have
\begin{equation*}
\mathbb{P}(\mathbf{Z}_{1}=\mathbf{z}_{1},\ldots ,\mathbf{Z}_{m}=\mathbf{z}%
_{m}|\mathbf{Z}_{0}=\mathbf{z};\mathcal{V}=v)=\mathcal{P}_{\mathbf{z},v}(%
\mathbf{Z}_{1}=\mathbf{z}_{1},\ldots ,\mathbf{Z}_{m}=\mathbf{z}_{m})\text{ \
\ a.s.}
\end{equation*}

In the sequel working in a random environment we use capital letters for the
variables or functions denoted by small letters in case of varying
environment. For instance, the (random) reproduction law of particles of the
$n$-th generation will be specified by the tuple $\mathbf{F}%
_{n}=(F_{n}^{(1)},\ldots,F_{n}^{(K)})$ of (random) probability generating
functions. Similarly, we denote by
\begin{equation}
\mathbf{M}_{n}:=(M_{n}(i,j))_{i,j=1}^{K}=\left( \frac{\partial F_{n}^{(i)}}{%
\partial s_{j}}(\mathbf{1})\right) _{i,j=1}^{K}  \label{matsr0}
\end{equation}%
the mean matrix corresponding to $\mathbf{F}_{n}$ and so on. Clearly, $%
\mathbf{M}_{n},n\geq 1,$ are i.i.d. matrices having the same distribution as
the matrix
\begin{equation}
\mathbf{M}=(M(i,j))_{i,j=1}^{K}=\left( \frac{\partial F^{(i)}}{\partial s_{j}%
}(\mathbf{1})\right) _{i,j=1}^{K}.  \label{matsr1}
\end{equation}

Now we formulate basic assumptions to be valid throughout the paper.

\textbf{Hypothesis 1}. We assume that the offspring generating functions of
our process are fractional-linear, that is, with probability 1 the
generating functions $\mathbf{F}$ and $\mathbf{F}_{n}$ have the form
\begin{equation}
\mathbf{F}(\mathbf{s})=\mathbf{1}-\dfrac{\mathbf{M}(\mathbf{1}-\mathbf{s})}{%
1+(\mathbf{w}\boldsymbol{,}\mathbf{1}-\mathbf{s})},\ \mathbf{F}_{n}(\mathbf{s%
})=\mathbf{1}-\dfrac{\mathbf{M}_{n}(\mathbf{1}-\mathbf{s})}{1+(\mathbf{w}%
_{n},\mathbf{1}-\mathbf{s})},\quad n\geq 1,  \label{FracLinear}
\end{equation}%
where $\boldsymbol{w}=\left( w_{1},\cdots ,w_{K}\right) $ and $\boldsymbol{w}%
_{j}=\left( w_{j1},\cdots ,w_{jK}\right) ,$ $j=1,\cdots ,K$ are independent $%
K$-dimensional nonnegative row random vectors of the environment with $%
\boldsymbol{w}\overset{d}{=}\boldsymbol{w}_{j}$.

\textbf{Hypothesis 2.} There exists a number $\alpha \in (0,1)$ such that
\begin{equation*}
\alpha \leq \dfrac{\max_{i,j}M(i,j)}{\min_{i,j}M(i,j)}\leq \frac{1}{\alpha }%
\quad
\end{equation*}%
with probability 1.

\textbf{Hypothesis 3. }There exists a nonrandom vector $\mathbf{v}=\left(
v_{1},\cdots ,v_{K}\right) $ with positive components such that, for $%
\mathbf{M}$ and all $\mathbf{M}_{n},n\geq 1$,
\begin{equation*}
\mathbf{vM}=\rho (\mathbf{M})\mathbf{v},\ \mathbf{vM}_{n}=\rho (\mathbf{M}%
_{n})\mathbf{v},
\end{equation*}%
that is, for matrices $\mathbf{M}$ and $\mathbf{M}_{n},n\geq 1,$ the $%
\mathbf{v}$ is the left eigenvector corresponding to their Perron roots.

Under this assumption, direct calculations with generating functions are
feasible, i.e. we can explicitly calculate the generating function of $%
\mathbf{Z}_{n}$, conditioned on the environment. Moreover, under this
condition a simple extension of the classification of the single type BPRE
given in~\cite{AGKV2005} to the MBPRE is possible.
Namely, setting \ $X:=\ln \rho ,$ $X_{n}:=\ln \rho _{n}$ and introducing the
so-called associated random walk $\mathcal{S}=\left\{ S_{n},n\geq 0\right\} $
as
\begin{equation*}
S_{0}:=0,\ S_{n}:=X_{1}+...+X_{n},
\end{equation*}%
we call a MBPRE satisfying Hypotheses 1 and 3 supercritical if $%
\lim\limits_{n\rightarrow \infty }S_{n}=\infty $ $a.s.,$ critical if $%
\limsup\limits_{n\rightarrow \infty }S_{n}=\infty $ $a.s.$ and $%
\liminf\limits_{n\rightarrow \infty }S_{n}=-\infty $ $a.s.,$ and subcritical
if $\lim\limits_{n\rightarrow \infty }S_{n}=-\infty $ $a.s.$

In particular, if $\mathbb{E}\left[ X\right] $ exists then a MBPRE $\left\{
\mathbf{Z}_{n},n\geq 0\right\} $ is supercritical, critical or subcritical,
if $\mathbb{E}\left[ X\right] >0$, $\mathbb{E}\left[ X\right] =0$ or $%
\mathbb{E}\left[ X\right] <0$, respectively.

The supercritical case admits additional classification based on the value
of $\mathbb{E}\left[ Xe^{-X}\right] $. \

We call a MBPRE meeting Hypothesis 3 strongly supercritical if $\mathbb{E}%
\left[ Xe^{-X}\right] >0$, intermediately supercritical if $\mathbb{E}\left[
Xe^{-X}\right] =0$ and weakly supercritical if $\mathbb{E}\left[ Xe^{-X}%
\right] <0$.

In the present paper, we mainly study the limiting behavior of the MBPRE in
the strongly and intermediately supercritical regimes and assume the
validity of the following

\textbf{Hypothesis 4}.
\begin{equation*}
\kappa :=\mathbb{E}\left[ e^{-X}\right] =\mathbb{E}\left[ \rho ^{-1}\right]
<\infty .
\end{equation*}

We use the parameter $\kappa $ to introduce a new measure $\mathbf{P}$
defined as follows: for any $n\in \mathbb{N}:=\{1,2,\ldots\}$ and any
bounded measurable function $\phi :(\mathfrak{P}(\mathbb{N}%
_{0}^{K}))^{n}\times (\mathbb{N}_{0}^{K})^{n+1}\rightarrow \mathbb{R}$,
\begin{equation*}
\mathbf{E}\left[ \phi \left( \mathbf{F}_{1},\ldots,\mathbf{F}_{n};\mathbf{Z}%
_{0},\ldots,\mathbf{Z}_{n}\right) \right] :=\kappa ^{-n}\mathbb{E}\left[
e^{-S_{n}}\phi \left( \mathbf{F}_{1},\ldots,\mathbf{F}_{n};\mathbf{Z}%
_{0},\ldots,\mathbf{Z}_{n}\right) \right] .
\end{equation*}

This change of measure is an important tool which will be used throughout
our paper.

Note that under the new measure the sequence $\left\{ \mathbf{Z}_{n},n\geq
0\right\} $ is still a $K$-type branching process in i.i.d. random
environment. Moreover, it is supercritical ($\mathbf{E}\left[ X\right] >0$)
for the initial strongly supercritical process, critical ($\mathbf{E}\left[ X%
\right] =0$) for the initial intermediately supercritical process and
subcritical ($\mathbf{E}\left[ X\right] <0$) for the initial weakly
supercritical process.

\section{Main results}

\subsection{The strongly supercritical case\label{Sub_strongly}}

In this subsection we assume that
\begin{equation}
0<\mathbf{E}\left[ X\right] =\mathbb{E}\left[ Xe^{-X}\right] <\infty ,
\label{strong}
\end{equation}%
there exists an $\varepsilon >0$ such that%
\begin{equation}
\mathbb{E}\left[ e^{-X}\left\vert X+\log \left\vert \mathbf{w}\right\vert
\right\vert ^{1+\varepsilon }\right] <\infty  \label{srrong01}
\end{equation}%
and, in addition, that
\begin{equation}
\mathbb{E}\left[ e^{-X}\log \left( \sum\limits_{i=1}^{K}\left( 1-\mathcal{P}%
_{\mathbf{e}_{i},v}\left( |\mathbf{Z}_{1}|=0\right) \right) \right) \right]
=\kappa \mathbf{E}\left[ \log \left( \sum\limits_{i=1}^{K}\left( 1-\mathcal{P%
}_{\mathbf{e}_{i},v}\left( |\mathbf{Z}_{1}|=0\right) \right) \right) \right]
>-\infty ,  \label{strong1}
\end{equation}%
which imply (see Theorem 9.10 in \cite{tanny}) that the process survives
under the measure $\mathbf{P}$ with a positive probability. Now we are ready
to formulate our first result concerning the limiting behavior of the MBPRE.

\begin{theorem}
\label{th1.1} If conditions (\ref{strong}) -- (\ref{strong1}) and Hypotheses
1 to 3 are valid then there exists a constant $\theta _{i}\in (0,\infty )$
such that, for each fixed $z=1,2,\ldots $%
\begin{equation}
\mathbb{P}_{\mathbf{e}_{i}}\left( \left\vert \mathbf{Z}_{n}\right\vert
=z\right) \thicksim \theta _{i}\cdot \kappa ^{n}  \label{LocalZ}
\end{equation}%
If, in addition,%
\begin{equation}
\mathbf{E}\left[ X\right] <\log \frac{1+\alpha ^{2}}{1-\alpha ^{2}},
\label{EgenValue0}
\end{equation}%
then, for each $\mathbf{z\in }\mathbb{N}_{0}^{K},\left\vert \mathbf{z}%
\right\vert >0$, there exists a constant $\Theta _{i}(\mathbf{z})>0$ such
that as $n\rightarrow \infty $,%
\begin{equation}
\mathbb{P}_{\mathbf{e}_{i}}\left( \mathbf{Z}_{n}=\mathbf{z}\right) \thicksim
\Theta _{i}(\mathbf{z})\cdot \kappa ^{n}.  \label{Vect}
\end{equation}
\end{theorem}

\begin{remark}
\label{Rem11} We see from (\ref{LocalZ}) that for every $c\in \mathbb{N}$
and $1\leq z\leq c$, as $n\rightarrow \infty $,
\begin{equation*}
\mathbb{P}_{\mathbf{e}_{i}}\left( |\mathbf{Z}_{n}|=z\Big|1\leq |\mathbf{Z}%
_{n}|\leq c\right) \rightarrow \frac{1}{c},
\end{equation*}%
i.e. the limiting distribution of $|\mathbf{Z}_{n}|$ is uniform on $%
\{1,\ldots ,c\}$.
\end{remark}

\begin{remark}
\label{Rem22}The explicit form of the constants $\Theta _{i}(\mathbf{z})$ is
given by formula (\ref{ExplicitTheta}) below.
\end{remark}

Before we state the next theorem, we formulate an important property
concerning the eigenvectors of the matrices
\begin{equation*}
\mathbf{M}_{k,n}:=\mathbf{M}_{k}\cdots \mathbf{M}_{n}=\left(
M_{k,n}(i,j)\right) _{i,j=1}^{K},\quad k=1,2,\ldots.
\end{equation*}%
We denote by $\mathbf{u}\left( \mathbf{M}_{1,n}\right) :=\left( u_{1}\left(
\mathbf{M}_{1,n}\right) ,,u_{K}\left( \mathbf{M}_{1,n}\right) \right) $ the
right positive eigenvector of $\mathbf{M}_{1,n}$ that corresponds to the
Perron root $\rho \left( \mathbf{M}_{1,n}\right) =\rho _{1}\cdots \rho _{n}$
and meets the scaling $\left( \mathbf{1},\mathbf{u}\left( \mathbf{M}%
_{1,n}\right) \right) =\left( \mathbf{v},\mathbf{u}\left( \mathbf{M}%
_{1,n}\right) \right) =1$.

\begin{lemma}
\label{u} (compare with points (i) and (ii) of Theorem 1 in \cite{Hen97})
Under Hypotheses 2 and 3, for all $i,j$ as $n\rightarrow \infty $%
\begin{equation*}
\frac{\mathbf{M}_{1,n}(i,j)}{\rho \left( \mathbf{M}_{1,n}\right) }%
-v_{j}u_{j}\left( \mathbf{M}_{1,n}\right) \rightarrow 0\ \mathbf{P}-a.s.,
\end{equation*}%
and there exists a random vector $\mathbf{u}=\left( u_{1},\ldots
,u_{K}\right) $ such that $\ \mathbf{u}\left( \mathbf{M}_{1,n}\right)
\rightarrow \mathbf{u}\ \ \mathbf{P}-a.s.$

In addition, $\left( \mathbf{v},\mathbf{u}\right) =1$.
\end{lemma}

With this lemma in hands, we can establish one more result concerning
conditional limiting properties of the process. For vectors $\mathbf{z}%
=\left( z_{1},\ldots ,z_{K}\right) $ and $\mathbf{v}=\left( v_{1},\ldots
,v_{K}\right) $ denote for brevity%
\begin{equation*}
C(\mathbf{z},\mathbf{v}):=\dfrac{\left\vert \mathbf{z}\right\vert !}{%
z_{1}!\cdots z_{K}!\left\vert \mathbf{v}\right\vert }\prod_{r=1}^{K}\left(
\frac{v_{r}}{|\mathbf{v}|}\right) ^{z_{r}},
\end{equation*}%
set%
\begin{equation*}
\mathcal{R}_{n}(i):=\mathcal{P}_{\mathbf{e}_{i}}(|\mathbf{Z}_{n}|=0)=%
\mathcal{P}_{\mathbf{e}_{i},v}(|\mathbf{Z}_{n}|=0),\ \mathcal{Q}_{n}(i):=%
\mathcal{P}_{\mathbf{e}_{i}}(|\mathbf{Z}_{n}|>0)=\mathcal{P}_{\mathbf{e}%
_{i},v}(|\mathbf{Z}_{n}|>0),
\end{equation*}%
and for the event \ $\{|\mathbf{Z}_{\infty }|>0\}:=\{|\mathbf{Z}%
_{n}|>0,\forall n\in \mathbb{N}\}$ let
\begin{equation*}
\mathcal{R}(i):=\mathcal{P}_{\mathbf{e}_{i}}(|\mathbf{Z}_{\infty }|=0),\
\mathcal{Q}(i):=1-\mathcal{R}(i)=\mathcal{P}_{\mathbf{e}_{i}}(|\mathbf{Z}%
_{\infty }|>0).
\end{equation*}

\begin{theorem}
\label{T_p}If conditions (\ref{strong}) -- (\ref{strong1}), (\ref{EgenValue0}) and Hypotheses
1 to 3 are valid then, for all $i,l\in \left\{ 1,\ldots ,K\right\} $ and $%
\mathbf{z}=(z_{1},\ldots ,z_{K})\neq \mathbf{0}$ and all $t\in (0,1)$,
\begin{equation*}
\lim\limits_{n\rightarrow \infty }\mathbb{P}_{\mathbf{e}_{i}}\left( \mathbf{%
Z}_{\lfloor nt\rfloor }=\mathbf{z}\Big|\mathbf{Z}_{n}=\mathbf{e}_{l}\right)
=p(\mathbf{z}),
\end{equation*}%
where $\mathbf{p}=(p(\mathbf{z}))_{\mathbf{z}\in \mathbb{N}_{0}^{K}}$ is a
proper probability distribution being the same for all $i$ and $l$ and
having components%
\begin{equation*}
p(\mathbf{z})=C(\mathbf{z},\mathbf{v})\mathbf{E}\left[ \sum_{r=1}^{K}\frac{%
z_{r}}{u_{r}}\frac{\mathcal{Q}^{2}(r)}{\mathcal{R}(r)}\cdot \prod_{j=1}^{K}%
\mathcal{R}^{z_{j}}(j)\right] .
\end{equation*}
\end{theorem}

\subsection{The intermediately supercritical case\label{Sub_interm}}

In this case, we assume, along with (\ref{strong1}) that
\begin{equation}
\mathbb{E}\left[ Xe^{-X}\right] =0.  \label{inter}
\end{equation}%
Then $\mathbf{E}\left[ X\right] =0$ and $\mathcal{S}$ is a recurrent random
walk under $\mathbf{P.}$

For our results, we impose some additional restrictions on the
characteristics of the MBPRE.

\textbf{Hypothesis 5.} The distribution of $X:=\log \rho $ under $\mathbf{P}$
belongs to the domain of attraction of a stable law with index $a\in (0,2]$.
It is a non-lattice distribution.

Define the random variable
\begin{equation*}
\vartheta :=1+\frac{1}{\rho ^{2}}\sum\limits_{\mathbf{z}\in \mathbb{N}%
^{K}}\sum\limits_{i=1}^{K}v_{i}\sum\limits_{j,k=1}^{K}z_{j}z_{k}F^{(i)}\left[
\mathbf{z}\right] .
\end{equation*}

\textbf{Hypothesis 6.} There is an $\varepsilon >0$ such that $ \mathbf{E}%
\left[ \log ^{a+\varepsilon }\vartheta \right] <\infty $, where the
parameter $a$ is from Hypothesis 5.

Now we are ready to formulate our next statement showing a different rate of
convergence to zero concerning the probabilities of rare events, compared with the
strongly supercritical case.

\begin{theorem}
\label{th2.1} If \ condition (\ref{inter}) and Hypotheses 1 to 6 are valid
then, for any $i\in \left\{ 1,2,...,K\right\} $ and $\mathbf{z}\in \mathbb{N}%
_{0}^{K},\left\vert \mathbf{z}\right\vert >0$, there exists a constant $%
\Delta _{i}(\mathbf{z})>0$ such that as $n\rightarrow \infty $,%
\begin{equation*}
\mathbb{P}_{\mathbf{e}_{i}}\left( \mathbf{Z}_{n}=\mathbf{z}\right) \thicksim
\Delta _{i}(\mathbf{z})\cdot \kappa ^{n}\mathbf{P}\left( \min\limits_{0\leq
k\leq n}S_{k}\geq 0\right) \sim \Delta _{i}(\mathbf{z})\cdot \kappa
^{n}n^{-(1-a)}l(n),
\end{equation*}%
where $l(n)$ is a function slowly varying at infinity.
\end{theorem}

The next theorem deals with the limiting distribution of $\{\mathbf{Z}%
_{n},n\geq 0\}$ when the values of $|\mathbf{Z}_{n}|$ are restricted to a
finite interval $[1,c]$.

\begin{theorem}
\label{T_uniformInterm}If condition (\ref{inter}) and Hypotheses 5 and 6 are
valid, then, for every $c\in \mathbb{N}$ and $1\leq z\leq c$,
\begin{equation*}
\lim_{n\rightarrow \infty }\mathbb{P}_{\mathbf{e}_{i}}\left( |\mathbf{Z}%
_{n}|=z\Big|1\leq |\mathbf{Z}_{n}|\leq c\right) =\frac{1}{c},
\end{equation*}%
i.e. the limiting distribution is uniform on $\{1,...,c\}$.
\end{theorem}

To state our next result, we need to introduce some notation related to the
associated random walk. Denote the time of the first minimum up to the
generation $n$ as
\begin{equation*}
\tau _{n}:=\min \left\{ 0\leq k\leq n:S_{k}=\min
\{S_{0},S_{1},...,S_{n}\}\right\} ,
\end{equation*}%
and the time of the first minimum between generations $k$ and $n$ by
\begin{equation*}
\tau _{k,n}:=\min \left\{ k\leq j\leq n:S_{j}=\min
\{S_{k},...,S_{n}\}\right\} .
\end{equation*}

\begin{theorem}
\label{T_cond_interm}If condition (\ref{inter}) and Hypotheses 1 to 6 are
valid then, for any $i,l\in \left\{ 1,...,K\right\} $ and $t\in (0,1)$
\begin{equation*}
\lim\limits_{n\rightarrow \infty }\mathbb{P}_{\mathbf{e}_{i}}\left( \mathbf{Z%
}_{\tau _{\lfloor nt\rfloor ,n}}=\mathbf{z}\Big|\mathbf{Z}_{n}=\mathbf{e}%
_{l}\right) =q(\mathbf{z}),
\end{equation*}%
where $\mathbf{q}=(q(\mathbf{z}))_{\mathbf{z}\in \mathbb{N}_{0}^{K}}$\ is a
probability distribution with
\begin{equation*}
q(\mathbf{z})=C(\mathbf{z},\mathbf{v})\mathbf{E}^{+}\left[
\sum\limits_{r=1}^{K}\dfrac{z_{r}}{u_{r}}\frac{\mathcal{Q}^{2}(r)}{\mathcal{R%
}(r)}\prod\limits_{j=1}^{K}\mathcal{R}^{z_{j}}(j)\right] .
\end{equation*}
\end{theorem}

\begin{remark}
\label{Rem2} The measure $\mathbf{P}^{+}$ corresponding to the expectation $%
\mathbf{E}^{+}$ will be defined in Section \ref{Sec_proof_interm}. Here we
just write the representation of $q(\mathbf{z})$.
\end{remark}

\begin{remark}
\label{Rem3} Note that the distributions $\mathbf{p}$ and $\mathbf{q}$ which
look similar are, in fact, essentially different: they are generated by
different measures.
\end{remark}

It follows from Theorem \ref{th2.1} that there exists a constant $\delta
_{i} $ such that, as $n\rightarrow \infty $
\begin{equation*}
\mathbb{P}_{\mathbf{e}_{i}}\left( |\mathbf{Z}_{n}|=1\right) \thicksim \delta
_{i}\cdot \kappa ^{n}n^{-(1-a)}l(n).
\end{equation*}

\subsection{Basic formulas\label{Sec_basic_form}}

Set $\mathbf{F}_{n,n}(\mathbf{s})=\mathbf{s,}$ and let, for $0\leq k<n$%
\begin{equation*}
\mathbf{F}_{k,n}(\mathbf{s})=\mathbf{F}_{k+1}\left( \mathbf{F}_{k+2}\cdots
\left( \mathbf{F}_{n}(\mathbf{s})\right) \right) =\left( F_{k,n}^{(1)}(%
\mathbf{s}),\ldots,F_{k,n}^{(K)}(\mathbf{s})\right) ,
\end{equation*}%
where $F_{k,n}^{(i)}(\mathbf{s})$ denote the probability generating function
for the number of particles in the $n$th generation given the process is
initiated at time $k$ by a single particle of type $i$. Recall that $\mathbf{%
F}_{n}(\mathbf{s})$ are all linear fractional generating functions. It is
easy to see by iterating that, for $n\geq 1$%
\begin{equation}
\mathbf{1}-\mathbf{F}_{0,n}\left( \mathbf{s}\right) =\dfrac{\mathbf{M}%
_{1}\left( \mathbf{1}-\mathbf{F}_{1,n}\left( \mathbf{s}\right) \right) }{1+(%
\mathbf{w}_{1},\mathbf{1}-\mathbf{F}_{1,n}\left( \mathbf{s}\right) )}=\dfrac{%
\mathbf{M}_{1,2}\left( \mathbf{1}-\mathbf{F}_{2,n}\left( \mathbf{s}\right)
\right) }{1+(\mathbf{w}_{1}\mathbf{M}_{2}+\mathbf{w}_{2},\mathbf{1}-\mathbf{F%
}_{2,n}\left( \mathbf{s}\right) )}=\ldots=\dfrac{\mathbf{M}_{1,n}\left(
\mathbf{1}-\mathbf{s}\right) }{1+(\mathbf{D}_{n},\mathbf{1}-\mathbf{s)}},
\label{2}
\end{equation}%
where
\begin{equation*}
\mathbf{D}_{n}:=\left( D_{n}(1),\ldots,D_{n}(K)\right) :=\mathbf{w}_{1}%
\mathbf{M}_{2,n}+\mathbf{w}_{2}\mathbf{M}_{3,n}+\cdots +\mathbf{w}_{n-1}%
\mathbf{M}_{n,n}+\mathbf{w}_{n}.
\end{equation*}%
Write
\begin{equation*}
\mathbf{M}_{1,n}(i)=\left( M_{1,n}(i,1),\ldots,M_{1,n}(i,K)\right) :=\mathbf{%
e}_{i}\mathbf{M}_{1,n}
\end{equation*}%
to denote the $i$th row of $\mathbf{M}_{1,n}$. It follows from (\ref{2}) that

\begin{eqnarray}
1-F_{0,n}^{(i)}(\mathbf{s}) &=&\dfrac{(\mathbf{e}_{i}\mathbf{M}_{1,n},%
\mathbf{1}-\mathbf{s)}}{1+(\mathbf{D}_{n},\mathbf{1}-\mathbf{s)}}=\dfrac{(%
\mathbf{M}_{1,n}(i),\mathbf{1}-\mathbf{s)}}{1+|\mathbf{D}_{n}|-(\mathbf{D}%
_{n},\mathbf{s)}}  \notag \\
&=&\dfrac{|\mathbf{M}_{1,n}(i)|}{1+|\mathbf{D}_{n}|}\sum\limits_{k=0}^{%
\infty }\left( \dfrac{(\mathbf{D}_{n},\mathbf{s)}}{1+|\mathbf{D}_{n}|}%
\right) ^{k}-\dfrac{(\mathbf{M}_{1,n}(i),\mathbf{s)}}{1+|\mathbf{D}_{n}|}%
\sum\limits_{k=0}^{\infty }\left( \dfrac{(\mathbf{D}_{n},\mathbf{s)}}{1+|%
\mathbf{D}_{n}|}\right) ^{k}.  \label{fn}
\end{eqnarray}

Substituting $\mathbf{s}=\mathbf{0}$, we get the explicit expressions for
survival and extinction probabilities:
\begin{equation}
\mathcal{Q}_{n}\left( i\right) =\mathcal{P}_{\mathbf{e}_{i},v}\left( |%
\mathbf{Z}_{n}|>0\right) =\dfrac{|\mathbf{M}_{1,n}(i)|}{1+|\mathbf{D}_{n}|}%
,\quad \mathcal{R}_{n}\left( i\right) =\mathcal{P}_{\mathbf{e}_{i},v}\left( |%
\mathbf{Z}_{n}|=0\right) =1-\dfrac{|\mathbf{M}_{1,n}(i)|}{1+|\mathbf{D}_{n}|}%
\ .  \label{4}
\end{equation}

For $z\geq 1,z\in \mathbb{N}$, by comparing the corresponding coefficient
for the total size of particles in the $n$th generation, regardless of the
types, we deduce that
\begin{equation}
\mathcal{P}_{\mathbf{e}_{i},v}\left( |\mathbf{Z}_{n}|=z\right) =\dfrac{|%
\mathbf{M}_{1,n}(i)|}{\left( 1+|\mathbf{D}_{n}|\right) ^{2}}\left( \dfrac{|%
\mathbf{D}_{n}|}{1+|\mathbf{D}_{n}|}\right) ^{z-1}.  \label{5}
\end{equation}

It is known \cite{kesten-spitzer} that for the product $\mathbf{M}%
_{k,n},k\leq n$ of the matrices satisfying Hypothesis 2, the relation
\begin{equation}
M_{k,n}(i,j)=\prod\limits_{r=k}^{n}\rho _{r}\cdot u_{i}(\mathbf{M}%
_{k,n})\cdot v_{j}\left( 1+\theta _{i,j}\beta ^{n-k}\right)  \label{Mnk}
\end{equation}%
is valid, where $\beta :=\left( 1-\alpha ^{2}\right) /\left( 1+\alpha
^{2}\right) \in (0,1)$, $|\theta _{i,j}|\leq C_{1}$, $C_{1}$ is a positive
constant. Further,
\begin{equation}
\left\vert \mathbf{M}_{1,n}(i)\right\vert
=\sum\limits_{j=1}^{K}M_{1,n}(i,j)\geq \frac{1}{\left\vert \mathbf{v}%
\right\vert }\sum\limits_{j=1}^{K}M_{1,n}(i,j)v_{j}=\frac{v_{i}}{\left\vert
\mathbf{v}\right\vert }e^{S_{n}}  \label{Bound_M}
\end{equation}%
and, similarly, for $v^{\ast }:=\min_{1\leq j\leq K}v_{j}$%
\begin{equation}
\left\vert \mathbf{M}_{1,n}(i)\right\vert \leq \frac{1}{v^{\ast }}%
\sum\limits_{j=1}^{K}M_{1,n}(i,j)v_{j}=\frac{v_{i}}{v^{\ast }}e^{S_{n}}.
\label{BoundM_above}
\end{equation}

Using the associated random walk $S_{n}=X_{1}+\cdots +X_{n}$, $n\geq 1$, $%
S_{0}=0$, with increments $X_{i}=\ln \rho _{i}$, $i\geq 1$, summing $%
M_{1,n}(i,j)$ over $j$ and applying (\ref{Mnk}), we see that
\begin{equation}
|\mathbf{M}_{1,n}(i)|=\sum\limits_{j=1}^{K}M_{1,n}(i,j)=|\mathbf{v}|u_{i}(%
\mathbf{M}_{1,n})\left( 1+O(\beta ^{n})\right) e^{S_{n}}.  \label{Mi}
\end{equation}%
Recalling the definition of $\mathbf{D}_{n}=\left( D_{n}(1),\ldots
,D_{n}(K)\right) $ and denoting
\begin{equation}
\eta _{k,n}:=\dfrac{|\mathbf{v}|}{\rho _{k}}\left( \mathbf{w}_{k},\mathbf{u}(%
\mathbf{M}_{k+1,n})\right)  \label{Def_eta_kn}
\end{equation}%
we conclude that \
\begin{eqnarray}
D_{n}(j) &=&\sum\limits_{k=1}^{n}\left( \mathbf{w}_{k},\mathbf{M}_{k+1,n}%
\mathbf{e}_{j}\right) =v_{j}\sum\limits_{k=1}^{n}\,\left( \mathbf{w}_{k},%
\mathbf{u}(\mathbf{M}_{k+1,n})\right) \prod\limits_{r=k+1}^{n}\rho _{r}\cdot
\left( 1+O(\beta ^{n-k})\right)  \notag \\
&=&\frac{v_{j}}{\left\vert \mathbf{v}\right\vert }\sum\limits_{k=1}^{n}\eta
_{k,n}e^{S_{n}-S_{k-1}}\left( 1+O(\beta ^{n-k})\right)  \label{D_indiv}
\end{eqnarray}%
and
\begin{equation}
|\mathbf{D}_{n}|=\sum\limits_{j=1}^{K}D_{n}(j)=\sum\limits_{k=1}^{n}\eta
_{k,n}e^{S_{n}-S_{k-1}}\left( 1+O(\beta ^{n-k})\right) .  \label{Dn}
\end{equation}

Combining (\ref{4}) and (\ref{5}) with the representation above and setting
\begin{equation}
H_{n}:=\dfrac{|\mathbf{D}_{n}|}{1+|\mathbf{D}_{n}|}=\dfrac{e^{-S_{n}}|%
\mathbf{D}_{n}|}{e^{-S_{n}}+e^{-S_{n}}|\mathbf{D}_{n}|},  \label{DefH}
\end{equation}%
we get for $z\geq 1$, $z\in \mathbb{N}$, $i=1,\cdots ,K,$ an important
formula which will be used many times in our proof,
\begin{equation}
\mathcal{P}_{\mathbf{e}_{i}}\left( |\mathbf{Z}_{n}|=z\right) =\mathcal{P}_{%
\mathbf{e}_{i},v}\left( |\mathbf{Z}_{n}|=z\right) =\dfrac{1}{|\mathbf{M}%
_{1,n}(i)|}\mathcal{Q}_{n}^{2}\left( i\right) H_{n}^{z-1}\quad \mathbf{P}%
-a.s.  \label{z}
\end{equation}

Note that $H_{n}$ is bounded by $1$ and if $S_{n}\rightarrow \infty \
\mathbf{P}-a.s.$ as $n\rightarrow \infty $, then also $H_{n}\rightarrow 1\
\mathbf{P}-a.s.$

\section{Proofs for the strongly supercritical case\label{SecProve_strong}}

In this section, we assume the validity of conditions (\ref{strong}) -- (\ref%
{strong1}). Let%
\begin{equation*}
\eta _{k}:=\dfrac{|\mathbf{v}|}{\rho _{k}}\left( \mathbf{w}_{k},\mathbf{u}%
^{(k+1)}\right) ,\quad k=1,2,\ldots,
\end{equation*}%
where $\mathbf{u}^{(k+1)}=(u_{1}^{(k+1)},\ldots ,u_{K}^{(k+1)})$ is the
limiting random vector, as $n\rightarrow \infty $ for $\mathbf{u}(\mathbf{M}%
_{k+1,n})$ in the sense of\ Lemma \ref{u}.

\begin{lemma}
\label{L_conergence}Given conditions (\ref{strong}) -- (\ref{strong1}) and
Hypotheses 2 to 4
\begin{equation}
\mathcal{G}:=\sum\limits_{k=0}^{\infty }\eta _{k+1}e^{-S_{k}}<\infty \quad
\mathbf{P}-a.s.  \label{FinG}
\end{equation}
\end{lemma}

\textbf{Proof}. \ By the strong law of large numbers $k^{-1}S_{k}\rightarrow
\mathbf{E}\left[ X\right] >0$ $\mathbf{P}$-a.s. as $k\rightarrow \infty $.
Hence, for any $\varepsilon >0$%
\begin{equation}
e^{-S_{k}}=O\left( e^{-k^{1/(1+\varepsilon /3)}}\right) \text{ as }%
k\rightarrow \infty \text{ \ \ }\mathbf{P}-a.s.  \label{Borel1}
\end{equation}
Further, in view of the inequality $\eta _{k}\leq |\mathbf{v}|\left\vert
\mathbf{w}_{k}\right\vert e^{-X_{k}}$ we have

\begin{eqnarray*}
\mathbf{P}\left( \eta _{k}>e^{k^{1/(1+\varepsilon /2)}}\right) &\leq &\left(
\frac{1}{k}\right) ^{\frac{1+\varepsilon }{1+\varepsilon /2}}\mathbf{E}\left[
\left( \log (|\mathbf{v}|\left\vert \mathbf{w}_{k}\right\vert
e^{-X_{k}})\right) ^{1+\varepsilon }\right] \\
&=&\left( \frac{1}{k}\right) ^{\frac{1+\varepsilon }{1+\varepsilon /2}}\frac{%
\mathbb{E}\left[ e^{-X}\left\vert X+\log |\mathbf{v}|\left\vert \mathbf{w}%
\right\vert \right\vert ^{1+\varepsilon }\right] }{\mathbb{E}\left[ e^{-X}%
\right] }.
\end{eqnarray*}%
Thus, by assumption (\ref{srrong01}) and Borel-Cantelli lemma%
\begin{equation*}
\eta _{k}=O\left( e^{k^{1/(1+\varepsilon /2)}}\right) \quad \mathbf{P}-a.s.
\end{equation*}%
This estimate combined with (\ref{Borel1}) justifies (\ref{FinG}).$\ \ \ \ \
\ \ \ \ \ \ \ \ \ \ \ \ \ \ \ \ \ \ \ \ \ \ \ \ \ \ \ \ \ \ \ \ \ \ \ \ \ \
\ \rule{2mm}{3mm}$ \hfill

\begin{lemma}
\label{L_eta(n)}Given conditions (\ref{strong}) -- (\ref{srrong01}) and
Hypotheses 2 to 4,
\begin{equation}
\lim_{n\rightarrow \infty }\sum\limits_{k=0}^{n}\eta _{k+1,n}e^{-S_{k}}=%
\mathcal{G}\quad \mathbf{P}-a.s.  \label{Eta_nk}
\end{equation}%
If, in addition, condition (\ref{EgenValue0}) is valid, then
\begin{equation}
\lim_{n\rightarrow \infty }e^{S_{n}}\sum\limits_{k=0}^{n}\eta
_{k+1,n}e^{-S_{k}}\beta ^{n-k}=0\quad \mathbf{P}-a.s.  \label{Negligible}
\end{equation}
\end{lemma}

\textbf{Proof}. \ Since, for each $k=1,2,\ldots$%
\begin{equation*}
\eta _{k,n}=\dfrac{|\mathbf{v}|}{\rho _{k}}\left( \mathbf{w}_{k},\mathbf{u}(%
\mathbf{M}_{k+1,n})\right) \rightarrow \dfrac{|\mathbf{v}|}{\rho _{k}}\left(
\mathbf{w}_{k},\mathbf{u}^{(k+1)}\right)
\end{equation*}%
as $n\rightarrow \infty $ in the sense of convergence established in Lemma %
\ref{u}, (\ref{Eta_nk}) follows from Lemma \ref{L_conergence}. To prove (\ref%
{Negligible}) it is necessary to observe that, as $n\rightarrow \infty $
\begin{equation*}
n^{-1}S_{n}+\log \beta \rightarrow \mathbf{E}\left[ X\right] -\log \frac{%
1+\alpha ^{2}}{1-\alpha ^{2}}<0\ \quad \mathbf{P}-a.s.
\end{equation*}%
and to repeat almost literally the proof of Lemma \ref{L_conergence}.\hfill
\rule{2mm}{3mm}\vspace{4mm}

\begin{lemma}
\label{L_ratioD}Given conditions (\ref{strong}) -- (\ref{srrong01}), and
Hypotheses 2 to 4,
\begin{equation*}
\lim_{n\rightarrow \infty }e^{-S_{n}}\mathbf{D}_{n}=\mathcal{G}\mathbf{%
v,\quad }\lim_{n\rightarrow \infty }\frac{\mathbf{D}_{n}}{1+|\mathbf{D}_{n}|}%
=\frac{\mathbf{v}}{\left\vert \mathbf{v}\right\vert }\quad \mathbf{P}-a.s.
\end{equation*}%
and%
\begin{equation*}
\lim_{n\rightarrow \infty }\frac{\mathbf{M}_{1,n}}{1+|\mathbf{D}_{n}|}=\frac{%
\mathbf{u}\otimes \mathbf{v}}{\left\vert \mathbf{v}\right\vert }\frac{1}{%
\mathcal{G}}\quad \mathbf{P}-a.s.,
\end{equation*}%
where $\mathbf{u}=(u_{1},\ldots ,u_{K})$ is a random vector given by Lemma %
\ref{u} and $\mathbf{u}\otimes \mathbf{v=}\left( u_{i}v_{j}\right)
_{i,j=1}^{K}$ is a $K\times K$ random matrix. In particular,%
\begin{equation}
\lim_{n\rightarrow \infty }\mathcal{Q}_{n}(i)=\lim_{n\rightarrow \infty }%
\frac{\left\vert \mathbf{M}_{1,n}(i)\right\vert }{1+|\mathbf{D}_{n}|}=\frac{%
u_{i}}{\mathcal{G}}\quad \mathbf{P}-a.s.  \label{AsympExtin}
\end{equation}
\end{lemma}

\textbf{Proof}. The first two statements are easy consequences of Lemma \ref%
{L_eta(n)}, equalities (\ref{D_indiv}) and (\ref{Dn}) and the evident fact
that $S_{n}\rightarrow \infty $ $\mathbf{P}-a.s.$ as $n\rightarrow \infty .~$%
To check the third and fourth statements it is sufficient to additionally
attract (\ref{Mnk}) with $k=1$. \rule{2mm}{3mm}\vspace{4mm}

\subsection{Proof of Theorem \protect\ref{th1.1}\label{Sec_provtheorem11}}

Recalling (\ref{z}) and the definition of the measure $\mathbf{P}$, we have
\begin{equation}
\mathbb{P}_{\mathbf{e}_{i}}(|\mathbf{Z}_{n}|=z)=\mathbb{E}\left[ \mathcal{P}%
_{\mathbf{e}_{i}}\left( |\mathbf{Z}_{n}|=z\right) \right] =\kappa ^{n}%
\mathbf{E}\left[ \dfrac{e^{S_{n}}}{\left\vert \mathbf{M}_{1,n}(i)\right\vert
}\mathcal{Q}_{n}^{2}\left( i\right) H_{n}^{z-1}\right] .  \label{Domin}
\end{equation}%
Since $\rho \left( \mathbf{M}_{1,n}\right) =e^{S_{n}}$, Lemma \ref{u} gives,
as $n\rightarrow \infty $%
\begin{equation}
\frac{\left\vert \mathbf{M}_{1,n}(i)\right\vert }{e^{S_{n}}}\rightarrow
u_{i}\left\vert \mathbf{v}\right\vert \quad \mathbf{P}-a.s.
\label{M_1nConverge}
\end{equation}%
Finally, $H_{n}\leq 1$ and $H_{n}\rightarrow 1$ $\mathbf{P}$ a.s., as $%
n\rightarrow \infty $. The estimates above and (\ref{Bound_M}) allow us to
apply the dominated convergence theorem to (\ref{Domin}) and to show by
Lemma \ref{L_conergence} that
\begin{eqnarray*}
\theta _{i}:= \lim\limits_{n\rightarrow \infty }\dfrac{\mathbb{P}_{\mathbf{e}%
_{i}}(|\mathbf{Z}_{n}|=z)}{\kappa ^{n}}&=&\mathbf{E}\left[
\lim\limits_{n\rightarrow \infty }\dfrac{e^{S_{n}}}{\left\vert \mathbf{M}%
_{1,n}(i)\right\vert }\lim\limits_{n\rightarrow \infty }\mathcal{Q}%
_{n}^{2}\left( i\right) \right] \\
&=&\mathbf{E}\left[ \dfrac{\mathcal{Q}^{2}\left( i\right) }{|\mathbf{v}|u_{i}%
}\right] =\frac{1}{|\mathbf{v}|}\mathbf{E}\left[ \dfrac{u_{i}}{\mathcal{G}%
^{2}}\right] \in (0,\infty ).
\end{eqnarray*}%
This completes the proof of the first part of the theorem.

To prove (\ref{Vect}) we fix $\mathbf{z}=(z_{1},\ldots ,z_{K})\in \mathbb{N}%
^{K},\left\vert \mathbf{z}\right\vert >0,$ and recalling the expression of $%
F_{0,n}^{(i)}(\mathbf{s})$ in (\ref{fn}) obtain, after some evident but
cumbersome transformations that
\begin{equation}
\quad \mathcal{P}_{\mathbf{e}_{i}}\left( \mathbf{Z}_{n}=\mathbf{z}\right) =%
\dfrac{\left\vert \mathbf{z}\right\vert !}{z_{1}!\cdots z_{K}!}%
\prod_{r=1}^{K}\left( \frac{D_{n}(r)}{1+|\mathbf{D}_{n}|}\right)
^{z_{r}}\left( \sum_{j=1}^{K}\frac{z_{j}M_{1,n}(i,j)}{\left\vert \mathbf{z}%
\right\vert D_{n}(j)}-\frac{|\mathbf{M}_{1,n}(i)|}{1+|\mathbf{D}_{n}|}%
\right) .  \label{SuperBasic}
\end{equation}%
We know by (\ref{Mnk}), (\ref{Mi}), (\ref{D_indiv}), (\ref{Dn}) and (\ref%
{Negligible}) that given (\ref{EgenValue0})%
\begin{equation}
\lim_{n\rightarrow \infty }e^{S_{n}}\left( \frac{M_{1,n}(i,j)}{|\mathbf{M}%
_{1,n}(i)|}\frac{|\mathbf{D}_{n}|}{D_{n}(j)}-1\right) =0\quad \mathbf{P}-a.s.
\label{Difference2}
\end{equation}%
Therefore,%
\begin{eqnarray}
&&\sum_{j=1}^{K}\frac{z_{j}M_{1,n}(i,j)}{\left\vert \mathbf{z}\right\vert
D_{n}(j)}-\frac{|\mathbf{M}_{1,n}(i)|}{1+|\mathbf{D}_{n}|}=|\mathbf{M}%
_{1,n}(i)|\left( \frac{1}{|\mathbf{D}_{n}|}\sum_{j=1}^{K}\frac{z_{j}}{%
\left\vert \mathbf{z}\right\vert }\frac{M_{1,n}(i,j)}{|\mathbf{M}_{1,n}(i)|}%
\frac{|\mathbf{D}_{n}|}{D_{n}(j)}-\frac{1}{1+|\mathbf{D}_{n}|}\right)  \notag
\\
&&  \notag \\
&&\qquad =|\mathbf{M}_{1,n}(i)|\left( \frac{1+o(e^{-S_{n}})}{|\mathbf{D}_{n}|%
}\sum_{j=1}^{K}\frac{z_{j}}{\left\vert \mathbf{z}\right\vert }-\frac{1}{1+|%
\mathbf{D}_{n}|}\right) =\frac{|\mathbf{M}_{1,n}(i)|}{|\mathbf{D}_{n}|(1+|%
\mathbf{D}_{n}|)}(1+\varepsilon _{i}(n)),  \label{Repres3}
\end{eqnarray}%
where $\varepsilon _{i}(n)\rightarrow 0$ $\mathbf{P}-$a.s. as $n\rightarrow
\infty .$ Using the equality%
\begin{equation*}
\mathbb{P}_{\mathbf{e}_{i}}\left( \mathbf{Z}_{n}=\mathbf{z}\right) =\mathbb{E%
}\left[ \mathcal{P}_{\mathbf{e}_{i}}\left( \mathbf{Z}_{n}=\mathbf{z}\right) %
\right] =\kappa ^{n}\mathbf{E}\left[ e^{S_{n}}\mathcal{P}_{\mathbf{e}%
_{i}}\left( \mathbf{Z}_{n}=\mathbf{z}\right) \right] ,
\end{equation*}%
applying the dominated convergence theorem (recall (\ref{4})) and Lemma \ref%
{L_ratioD} we obtain%
\begin{eqnarray}
\lim\limits_{n\rightarrow \infty }\dfrac{\mathbb{P}_{\mathbf{e}_{i}}(%
\mathbf{Z}_{n}=\mathbf{z})}{\kappa ^{n}}&=&\dfrac{\left\vert \mathbf{z}%
\right\vert !}{z_{1}!\cdots z_{K}!}\mathbf{E}\left[ \lim_{n\rightarrow
\infty }\prod_{r=1}^{K}\left( \frac{D_{n}(r)}{1+|\mathbf{D}_{n}|}\right)
^{z_{r}}\frac{|\mathbf{M}_{1,n}(i)|}{1+|\mathbf{D}_{n}|}\frac{e^{S_{n}}}{|%
\mathbf{D}_{n}|}\right]  \notag \\
&=&\dfrac{\left\vert \mathbf{z}\right\vert !}{z_{1}!\cdots z_{K}!}%
\prod_{r=1}^{K}\left( \frac{v_{r}}{|\mathbf{v}|}\right) ^{z_{r}}\mathbf{E}%
\left[ \frac{u_{i}}{|\mathbf{v}|\mathcal{G}^{2}}\right] \notag\\
&=&\dfrac{\left\vert
\mathbf{z}\right\vert !}{z_{1}!\cdots z_{K}!|\mathbf{v}|}\prod_{r=1}^{K}%
\left( \frac{v_{r}}{|\mathbf{v}|}\right) ^{z_{r}}\mathbf{E}\left[ \frac{%
\mathcal{Q}^{2}(i)}{u_{i}}\right] =:\Theta _{i}(\mathbf{z}).
\label{ExplicitTheta}
\end{eqnarray}
The proof of Theorem \ref{th1.1} is completed.\hfill \rule{2mm}{3mm}\vspace{%
4mm}

\subsection{Proof of Theorem \protect\ref{T_p}\label{SubT_p}}

Fix $\mathbf{z}\in \mathbb{N}^{K}$ , $\mathbf{z}\neq \mathbf{0},$ and $t\in
(0,1)$. By the Markov property of the process and the independence of the
environment we have
\begin{align}
\mathbb{P}_{\mathbf{e}_{i}}\left( \mathbf{Z}_{\lfloor nt\rfloor }=\mathbf{z},%
\mathbf{Z}_{n}=\mathbf{e}_{l}\right) & =\mathbb{E}\left[ \mathcal{P}_{%
\mathbf{e}_{i}}\left( \mathbf{Z}_{\lfloor nt\rfloor }=\mathbf{z},\mathbf{Z}%
_{n}=\mathbf{e}_{l}\right) \right]  \notag \\
& =\sum\limits_{\left\vert \mathbf{z}\right\vert =z}\mathbb{P}_{\mathbf{e}%
_{i}}\left( \mathbf{Z}_{\lfloor nt\rfloor }=\mathbf{z}\right) \cdot \mathbb{P%
}_{\mathbf{z}}\left( \mathbf{Z}_{n-\lfloor nt\rfloor }=\mathbf{e}_{l}\right)
,  \notag
\end{align}%
where $\mathbb{P}_{\mathbf{z}}(\cdot )=\mathbb{P}_{\left( z_{1},\cdots
,z_{K}\right) }(\cdot )$ denotes the distribution of the process with the
initial state $\mathbf{Z}_{0}=\mathbf{z}=\left( z_{1},\cdots ,z_{K}\right) $.

From Theorem \ref{strong}, we know that, as $n\rightarrow \infty $
\begin{equation}
\mathbb{P}_{\mathbf{e}_{i}}\left( \mathbf{Z}_{\lfloor nt\rfloor }=\mathbf{z}%
\right) \sim \kappa ^{\lfloor nt\rfloor }\cdot C(\mathbf{z},\mathbf{v})%
\mathbf{E}\left[ \frac{\mathcal{Q}^{2}(i)}{u_{i}}\right] .  \label{nt}
\end{equation}%
Besides, in view of (\ref{SuperBasic}) -- (\ref{Repres3})
\begin{eqnarray*}
\mathcal{P}_{\mathbf{e}_{j}}\left( \mathbf{Z}_{n-\lfloor nt\rfloor }=\mathbf{%
e}_{l}\right) &=&\frac{D_{n-\lfloor nt\rfloor }(l)}{\left\vert \mathbf{D}%
_{n-\lfloor nt\rfloor }\right\vert }\frac{\left\vert \mathbf{M}_{1,n-\lfloor
nt\rfloor }(i)\right\vert }{\left( 1+\left\vert \mathbf{D}_{n-\lfloor
nt\rfloor }\right\vert \right) ^{2}}\left( 1+\varepsilon _{i}(n-\lfloor
nt\rfloor )\right) \\
&=&\frac{D_{n-\lfloor nt\rfloor }(l)}{\left\vert \mathbf{D}_{n-\lfloor
nt\rfloor }\right\vert }\frac{1}{\left\vert \mathbf{M}_{1,n-\lfloor
nt\rfloor }(i)\right\vert }\mathcal{Q}_{n-\lfloor nt\rfloor }^{2}(j)\left(
1+\varepsilon _{i}(n-\lfloor nt\rfloor )\right) .
\end{eqnarray*}

Starting from the state $\mathbf{Z}_{0}=\mathbf{z}$, the event $\{\mathbf{Z}%
_{n-\lfloor nt\rfloor }=\mathbf{e}_{l}\}$ means that only a single
subpopulation steamed from the initial $\left\vert \mathbf{z}\right\vert $
particles survives, while the other $\left\vert \mathbf{z}\right\vert -1$
subpopulations extinct before time $n-\lfloor nt\rfloor $. So we get by (\ref%
{z}) as $n\rightarrow \infty $
\begin{eqnarray}
\mathcal{P}_{\mathbf{z}}\left( \mathbf{Z}_{n-\lfloor nt\rfloor }=\mathbf{e}%
_{l}\right) &=&\sum_{j=1}^{K}\left[ z_{j}\mathcal{R}_{n-\lfloor nt\rfloor
}^{z_{j}-1}(j)\mathcal{P}_{\mathbf{e}_{j}}\left( \mathbf{Z}_{n-\lfloor
nt\rfloor }=\mathbf{e}_{l}\right) \prod_{k\neq j}\mathcal{R}_{n-\lfloor
nt\rfloor }^{z_{k}}(k)\right]  \notag \\
&\sim &\frac{D_{n-\lfloor nt\rfloor }(l)}{\left\vert \mathbf{D}_{n-\lfloor
nt\rfloor }\right\vert }\sum_{j=1}^{K}\Bigg[\frac{z_{j}}{\left\vert \mathbf{M%
}_{1,n-\lfloor nt\rfloor }(i)\right\vert }\frac{\mathcal{Q}_{n-\lfloor
nt\rfloor }^{2}(j)}{\mathcal{R}_{n-\lfloor nt\rfloor }(j)}\cdot
\prod_{k=1}^{K}\mathcal{R}_{n-\lfloor nt\rfloor }^{z_{k}}(k)\Bigg].
\label{Z_nt}
\end{eqnarray}%
Using the change of measure, we obtain
\begin{equation}
\frac{\mathbb{P}_{\mathbf{z}}\left( \mathbf{Z}_{n-\lfloor nt\rfloor }=%
\mathbf{e}_{l}\right) }{\kappa ^{n-\lfloor nt\rfloor }}\sim \mathbf{E}\Bigg[%
\frac{D_{n-\lfloor nt\rfloor }(l)}{\left\vert \mathbf{D}_{n-\lfloor
nt\rfloor }\right\vert }\sum_{j=1}^{K}\frac{z_{j}e^{S_{n-\lfloor nt\rfloor }}%
}{\left\vert \mathbf{M}_{1,n-\lfloor nt\rfloor }(i)\right\vert }\frac{%
\mathcal{Q}_{n-\lfloor nt\rfloor }^{2}(j)}{\mathcal{R}_{n-\lfloor nt\rfloor
}(j)}\prod_{k=1}^{K}\mathcal{R}_{n-\lfloor nt\rfloor }^{z_{k}}(k)\Bigg].
\notag
\end{equation}%
By the dominated convergence theorem, Lemma \ref{L_ratioD} and the
inequality $\mathbf{P}\left( |\mathbf{Z}_{\infty }|>0\right) >0$ we conclude
that
\begin{equation*}
\quad \frac{\mathbb{P}_{\mathbf{z}}\left( \mathbf{Z}_{n-\lfloor nt\rfloor }=%
\mathbf{e}_{l}\right) }{\kappa ^{n-\lfloor nt\rfloor }}\sim \frac{v_{l}}{%
\left\vert \mathbf{v}\right\vert ^{2}}\mathbf{E}\left[ \sum_{j=1}^{K}\frac{%
z_{j}}{u_{j}}\frac{\mathcal{Q}^{2}(j)}{\mathcal{R}(j)}\cdot \prod_{k=1}^{K}%
\mathcal{R}^{z_{k}}(k)\right] >0,
\end{equation*}%
as $n\rightarrow \infty $. By (\ref{ExplicitTheta}), we also know that
\begin{equation}
\frac{\mathbb{P}_{\mathbf{e}_{i}}\left( \mathbf{Z}_{n}=\mathbf{e}_{l}\right)
}{\kappa ^{n}}\sim \frac{v_{l}}{\left\vert \mathbf{v}\right\vert ^{2}}%
\mathbf{E}\left[ \frac{\mathcal{Q}^{2}(i)}{u_{i}}\right] ,  \label{n}
\end{equation}%
as $n\rightarrow \infty $. Combining (\ref{nt}) and (\ref{n}) we get the
following relations:%
\begin{eqnarray*}
\mathbb{P}_{\mathbf{e}_{i}}\left( \mathbf{Z}_{\lfloor nt\rfloor }=\mathbf{z}%
\Big|\mathbf{Z}_{n}=\mathbf{e}_{l}\right) &=&\dfrac{\mathbb{P}_{\mathbf{e}%
_{i}}\left( \mathbf{Z}_{\lfloor nt\rfloor }=\mathbf{z}\right) \cdot \mathbb{P%
}_{\mathbf{z}}\left( \mathbf{Z}_{n-\lfloor nt\rfloor }=\mathbf{e}_{l}\right)
}{\mathbb{P}_{\mathbf{e}_{i}}\left( \mathbf{Z}_{n}=\mathbf{e}_{l}\right) } \\
&\sim &C(\mathbf{z},\mathbf{v})\mathbf{E}\left[ \sum_{j=1}^{K}\frac{z_{j}}{%
u_{j}}\frac{\mathcal{Q}^{2}(j)}{\mathcal{R}(j)}\cdot \prod_{k=1}^{K}\mathcal{%
R}^{z_{k}}(k)\right] =:p(\mathbf{z}).
\end{eqnarray*}%
Note that $p(\mathbf{z})$ does not depend on $i,l$ and $t$.

Let us check that $\mathbf{p}:=\left( p(\mathbf{z})\right) _{\mathbf{z}\in
\mathbb{N}_{0}^{K}}$ is a proper distribution. It is not difficult to see
that
\begin{align}
& \quad \sum\limits_{\left\vert \mathbf{z}\right\vert =z}C(\mathbf{z},%
\mathbf{v})\mathbf{E}\left[ \sum_{j=1}^{K}\frac{z_{j}}{u_{j}}\frac{\mathcal{Q%
}^{2}(j)}{\mathcal{R}(j)}\cdot \prod_{k=1}^{K}\mathcal{R}^{z_{k}}(k)\right]
\notag \\
& =\mathbf{E}\left[ \sum\limits_{\left\vert \mathbf{z}\right\vert
=z}\sum_{j=1}^{K}C(\mathbf{z},\mathbf{v})\frac{z_{j}}{u_{j}}\mathcal{Q}%
^{2}(j)\mathcal{R}^{z_{j}-1}(j)\prod_{k\neq j}\mathcal{R}^{z_{k}}(k)\right]
\notag \\
& =\mathbf{E}\Bigg[\sum_{j=1}^{K}\frac{\mathcal{Q}^{2}(j)}{u_{j}}%
\sum\limits_{\left\vert \mathbf{z}\right\vert =z}\dfrac{z\cdot (z-1)!v_{j}}{%
z_{1}!\cdots (z_{j}-1)!\cdots z_{K}!|\mathbf{v}|}\left( \frac{v_{j}}{|%
\mathbf{v}|}\mathcal{R}(j)\right) ^{z_{j}-1}\prod_{k\neq j}\left( \frac{v_{k}%
}{|\mathbf{v}|}\mathcal{R}(k)\right) ^{z_{k}}\Bigg]  \notag \\
& =\mathbf{E}\left[ \sum_{j=1}^{K}\frac{\mathcal{Q}^{2}(j)}{u_{j}}\cdot z%
\frac{v_{j}}{|\mathbf{v}|}\left( \frac{v_{1}}{|\mathbf{v}|}\mathcal{R}%
(1)+\cdots +\frac{v_{K}}{|\mathbf{v}|}\mathcal{R}(K)\right) ^{z-1}\right]
\notag \\
& =z\mathbf{E}\left[ \left( \sum_{i=1}^{K}\frac{v_{i}}{\left\vert \mathbf{v}%
\right\vert }\mathcal{R}(i)\right) ^{z-1}\sum_{j=1}^{K}\dfrac{v_{j}}{|%
\mathbf{v}|u_{j}}\mathcal{Q}^{2}(j)\right]  \notag
\end{align}%
implying
\begin{equation*}
\lim\limits_{n\rightarrow \infty }\mathbb{P}_{\mathbf{e}_{i}}\left( |\mathbf{%
Z}_{\lfloor nt\rfloor }|=z\Big|\mathbf{Z}_{n}=\mathbf{e}_{l}\right) =\dfrac{z%
}{|\mathbf{v}|}\mathbf{E}\left[ \left( \sum_{i=1}^{K}\frac{v_{i}}{|\mathbf{v}%
|}\mathcal{R}(i)\right) ^{z-1}\sum_{j=1}^{K}\dfrac{v_{j}}{u_{j}}\mathcal{Q}%
^{2}(j)\right] .
\end{equation*}

Using, for $x\in \lbrack 0,1)$, the equality $\sum\limits_{z=1}^{\infty
}zx^{z-1}=\dfrac{1}{(1-x)^{2}}$ and observing that%
\begin{equation*}
\sum_{k=1}^{K}\frac{v_{k}}{\left\vert \mathbf{v}\right\vert }\mathcal{R}%
(k)=\sum_{k=1}^{K}\frac{v_{k}}{\left\vert \mathbf{v}\right\vert }\left( 1-%
\mathcal{Q}(k)\right) =\sum_{k=1}^{K}\frac{v_{k}}{\left\vert \mathbf{v}%
\right\vert }\left( 1-\frac{u_{k}}{\mathcal{G}}\right) =1-\frac{1}{%
\left\vert \mathbf{v}\right\vert \mathcal{G}}
\end{equation*}%
and
\begin{equation*}
\sum_{j=1}^{K}\dfrac{v_{j}}{|\mathbf{v}|u_{j}}\mathcal{Q}^{2}(j)=\frac{1}{|%
\mathbf{v}|\mathcal{G}^{2}}\sum_{j=1}^{K}v_{j}u_{j}=\frac{1}{|\mathbf{v}|%
\mathcal{G}^{2}}
\end{equation*}%
we see that
\begin{equation}
\sum_{\mathbf{z}\in \mathbb{N}_{0}^{K}}^{\infty }p(\mathbf{z}%
)=\sum_{z=1}^{\infty }\mathbf{E}\Bigg[\dfrac{z}{|\mathbf{v}|^{2}}\dfrac{1}{%
\mathcal{G}^{2}}\left( 1-\frac{1}{\left\vert \mathbf{v}\right\vert \mathcal{G%
}}\right) ^{z-1}\Bigg]=1.  \notag
\end{equation}%
This completes the proof of Theorem \ref{T_p}.\hfill \rule{2mm}{3mm}\vspace{%
4mm}

\section{Proofs of the intermediately supercritical case\label%
{Sec_proof_interm}}

In this section, we assume that $\mathbb{E}\left[ Xe^{-X}\right] =0$, and
thus $\mathbf{E}\left[ X\right] =0$.

\subsection{Some properties of random walks}

In the proofs concerning the intermediately supercritical case, we analyze
the trajectory of the environment in more detail using mainly the tools of
the associated random walk $\mathcal{S}$. To describe them, we first
introduce some notation and list needed properties of $\mathcal{S}$.

For $0\leq k\leq n$, we define
\begin{equation*}
L_{k,n}:=\min\limits_{0\leq j\leq n-k}\left\{ S_{k+j}-S_{k}\right\} ,\quad
L_{n}:=L_{0,n}=\min \left\{ S_{0},\ldots ,S_{n}\right\} ,\quad M_{n}:=\max
\left\{ S_{0},\ldots ,S_{n}\right\} .
\end{equation*}%
Recall the definition of $\tau _{n}$, $\tau _{n}=\min \left\{ 0\leq k\leq
n:S_{k}=L_{n}\right\} $.

Introduce the renewal function $V(x)$ specified by the relations
\begin{equation*}
V(x):=1+\sum\limits_{k=1}^{\infty }\mathbf{P}(-S_{k}\leq x,M_{k}<0)\text{ if
}~x\geq 0,
\end{equation*}%
and zero otherwise. Since $\mathbf{E}\left[ X\right] =0$, it follows that
(see, for instance, \cite{AGKV2005})
\begin{equation}
\mathbf{E}\left[ V(x+X);x+X\geq 0\right] =V(x),\quad x\geq 0.
\label{martingale}
\end{equation}%
Below we describe properties of random walks conditioned to never hit the
strictly negative half line. To do this in a formal way we need to specify
an auxiliary measure $\mathbf{P^{+}}$ as follows.

Let $\{\mathcal{F}_{n},n\in \mathbb{N}_{0}\}$ be a sequence of $\sigma $%
-algebras defined by
\begin{equation*}
\mathcal{F}_{0}:=\sigma \left( \mathbf{Z}_{0}\right) ,\qquad \mathcal{F}%
_{n}:=\sigma \left( \mathbf{F}_{1},\ldots ,\mathbf{F}_{n},\mathbf{Z}%
_{0},\ldots ,\mathbf{Z}_{n}\right) ,\quad n\geq 1.
\end{equation*}%
Using (\ref{martingale}) we introduce a measure $\mathbf{P^{+}}$ by setting
for any bounded $\mathcal{F}_{n}$-measurable random variable $Y_{n}$
\begin{equation*}
\mathbf{E^{+}}Y_{n}:=\mathbf{E}\left[ Y_{n}V(S_{n});L_{n}\geq 0\right] ,
\end{equation*}%
where $\mathbf{E^{+}}$ is the expectation corresponding to the measure $%
\mathbf{P^{+}}$. The reader may find in \cite{AGKV2005} more details showing
that $\mathbf{P^{+}}$ is well defined.

For a probability measure $\mu \in \mathbb{N}_{0}^{K}$, we use the symbol $%
\mathbf{Z}_{n}^{\mu }$, $n\geq 0$, to denote a $K$-type MBPRE if $\mathbf{Z}%
_{0}^{\mu }$ is distributed as $\mu $, and use the symbols $\mathbb{P}_{\mu
} $, $\mathbb{E}_{\mu }$ to denote the corresponding probability law and the
expectation. With this notation in hands we now are ready to formulate
several important statements.

\begin{lemma}
\textrm{{\textbf{\cite{AGKV2005}}}}\label{dya1} Under (\ref{inter}), let $%
\mu $ be a probability measure on $\mathbb{N}_{0}^{K}$ which is not
concentrated at $\mathbf{0}$. Let $Y_{n}$, $n\in \mathbb{N}_{0}$, be a
uniformly bounded sequence of real-valued random variables adapted to the
filtration $\{\mathcal{F}_{n}\}$ which converges $\mathbf{P}_{\mu }^{+}-a.s.$
to some random variable $Y_{\infty }$. Then as $n\rightarrow \infty $,
\begin{equation*}
\mathbf{E}_{\mu }\left[ Y_{n}\Big|L_{n}\geq 0\right] \rightarrow \mathbf{E}%
_{\mu }^{+}\left[ Y_{\infty }\right] .
\end{equation*}
\end{lemma}

\begin{lemma}
\textrm{{\textbf{\cite{AGKV2005}}}}\label{dya2} Under (\ref{inter}), suppose
that $V_{1},V_{2},\ldots $ is a uniformly bounded sequence of real-valued
random variables which for every $k\geq 0$ satisfy the relation
\begin{equation*}
\mathbf{E}\left[ V_{n};|\mathbf{Z}_{k}|>0,L_{k,n}\geq 0\Big|\mathcal{F}_{k}%
\right] =\mathbf{P}\left( L_{n}\geq 0\right) \left( V_{k,\infty
}+o(1)\right) \quad \mathbf{P}-a.s.
\end{equation*}%
for some sequence of random variables $V_{1,\infty },V_{2,\infty },\ldots .$
Then
\begin{equation*}
\mathbf{E}\left[ V_{n};|\mathbf{Z}_{\tau _{n}}|>0\right] =\mathbf{P}%
(L_{n}\geq 0)\left( \sum\limits_{k=0}^{\infty }\mathbf{E}\left[ V_{k,\infty
};\tau _{k}=k\right] +o(1)\right) ,
\end{equation*}%
where the series in the right-hand side is absolutely convergent.
\end{lemma}

\begin{lemma}
\label{l7} For every $\varepsilon>0 $ and $z\in\mathbb{N} $ there exists an $%
m=m(\varepsilon)\in\mathbb{N} $ such that
\begin{equation*}
\mathbb{P}_{\mathbf{e}_{i}}\left( \vert\mathbf{Z}_{n}\vert=z,
\tau_{n}>m\right)<\varepsilon\kappa^{n}\mathbf{P}(L_{n}\geq 0),\qquad
i=1,\ldots,K.
\end{equation*}
\end{lemma}

\textbf{Proof} By (\ref{z}) and (\ref{Bound_M}), we have
\begin{eqnarray*}
&&\mathbb{P}_{\mathbf{e}_{i}}\left( |\mathbf{Z}_{n}|=z,\tau _{n}>m\right) =%
\mathbb{E}\left[ \mathcal{P}_{\mathbf{e}_{i}}\left( |\mathbf{Z}%
_{n}|=z\right) ;\tau _{n}>m\right] \\
&&\qquad =\mathbb{E}\left[ \dfrac{1}{\left\vert \mathbf{M}%
_{1,n}(i)\right\vert }\mathcal{Q}_{n}^{2}\left( i\right) H_{n}^{z-1};\tau
_{n}>m\right] \leq \frac{\left\vert \mathbf{v}\right\vert }{v_{i}}\mathbb{E}%
\left[ e^{-S_{n}}\mathcal{Q}_{n}^{2}\left( i\right) ;\tau _{n}>m\right] .
\end{eqnarray*}%
We know from (\ref{4}) and (\ref{BoundM_above}) that
\begin{equation*}
\mathcal{Q}_{n}(i)\leq \min_{0\leq k\leq n}\mathcal{Q}_{k}(i)\leq
\min_{0\leq k\leq n}\left\vert \mathbf{M}_{1,k}(i)\right\vert \leq \frac{%
v_{i}}{v^{\ast }}\min_{0\leq k\leq n}e^{S_{k}}=\frac{v_{i}}{v^{\ast }}%
e^{L_{n}}\qquad \mathbb{P}-a.s.
\end{equation*}%
Recall the definition of the measure $\mathbf{P}$ and use the fact that $%
\{\tau _{n}=k\}=\{\tau _{k}=k,L_{k,n}\geq 0\}$,
\begin{eqnarray*}
&&\mathbb{P}_{\mathbf{e}_{i}}\left( |\mathbf{Z}_{n}|=z,\tau _{n}>m\right)
\leq \frac{\left\vert \mathbf{v}\right\vert }{v^{\ast }}\sum_{k=m+1}^{n}%
\mathbb{E}\left[ e^{-S_{n}+2L_{n}};\tau _{k}=k,L_{k,n}\geq 0\right] \\
&=&\frac{\left\vert \mathbf{v}\right\vert }{v^{\ast }}\kappa
^{n}\sum_{k=m+1}^{n}\mathbf{E}\left[ e^{2L_{n}};\tau _{k}=k,L_{k,n}\geq 0%
\right] =\frac{\left\vert \mathbf{v}\right\vert }{v^{\ast }}\kappa
^{n}\sum_{k=m+1}^{n}\mathbf{E}\left[ e^{2S_{k}};\tau _{k}=k\right] \mathbf{P}%
\left( L_{n-k}\geq 0\right) .
\end{eqnarray*}%
From Lemma 2.2 in \cite{AGKV2005} by letting $%
u(x)=e^{-2x}$, we deduce that for every $\varepsilon >0$ and $m\in \mathbb{N}
$ big enough,
\begin{equation*}
\sum_{k=m+1}^{n}\mathbf{E}\left[ e^{2S_{k}};\tau _{k}=k\right] \mathbf{P}%
\left( L_{n-k}\geq 0\right) <\varepsilon \mathbf{P}(L_{n}\geq 0).
\end{equation*}%
Hence the desired statement follows.\hfill \rule{2mm}{3mm}\vspace{4mm}

\begin{remark}
\label{Rem1}It follows from the proof of the lemma, that for every $%
\varepsilon >0$ and sufficiently large $r$
\begin{equation*}
\limsup_{n\rightarrow \infty }\mathbf{P}(L_{n}\geq 0)^{-1}\sum_{k=r+1}^{n}%
\mathbf{E}\left[ \mathcal{Q}_{k}^{2}(i);\tau _{n}=k\right] <\varepsilon .
\end{equation*}
\end{remark}

\begin{lemma}
\label{Zj} For every $k$ and $\mathbf{z}=(z_{1},\ldots ,z_{K})\in \mathbb{N}%
_{0}^{K}$,
\begin{equation*}
T(\mathbf{z}):=\lim\limits_{n\rightarrow \infty }\mathbf{E}\left[ e^{S_{n}}%
\mathcal{P}_{\mathbf{z}}\left( |\mathbf{Z}_{n}|=k\right) \Big|L_{n}\geq 0%
\right] =\mathbf{E}^{+}\left[ \sum\limits_{i=1}^{K}\dfrac{z_{i}}{|\mathbf{v}%
|u_{i}}\frac{\mathcal{Q}^{2}(i)}{\mathcal{R}(i)}\prod\limits_{j=1}^{K}%
\mathcal{R}^{z_{j}}(j)\right] ,
\end{equation*}

where the limit does not depend on $k$.
\end{lemma}

\textbf{Proof} We will complete the proof by two steps and mainly use the
method of induction.

\noindent \textit{step 1, induction in single type. }

First we prove that for each $z_{i}\geq 1$
\begin{equation}
e^{S_{n}}\mathcal{P}_{z_{i}\mathbf{e}_{i}}\left( |\mathbf{Z}_{n}|=k\right)
\overset{n\rightarrow \infty }{\longrightarrow }\dfrac{z_{i}}{|\mathbf{v}%
|u_{i}}\mathcal{Q}^{2}(i)\mathcal{R}^{z_{i}-1}(i)\qquad \mathbf{P}^{+}-a.s.
\label{zi}
\end{equation}%
We know that under $\mathbf{P}^{+}$, $S_{n}\rightarrow \infty $ a.s. This
and (\ref{DefH}) show that $H_{n}\rightarrow 1$ $\mathbf{P}^{+}$-a.s.

Recalling now (\ref{z}) and (\ref{Mi}) gives
\begin{equation*}
\lim\limits_{n\rightarrow \infty }e^{S_{n}}\mathcal{P}_{\mathbf{e}%
_{i}}\left( |\mathbf{Z}_{n}|=k\right) =\lim\limits_{n\rightarrow \infty }%
\dfrac{1}{|\mathbf{v}|u_{i}}\mathcal{Q}_{n}^{2}\left( i\right) H_{n}^{k-1}=%
\dfrac{1}{|\mathbf{v}|u_{i}}\mathcal{Q}^{2}(i)\qquad \mathbf{P}^{+}-a.s.,
\end{equation*}%
proving the desired statement for $z_{i}=1$. Assume that (\ref{zi}) is valid
for some $z_{i}\in \mathbb{N}$. For $z_{i}+1$, the process starts from $%
z_{i}+1$ particles of type $i$, and the size of the $n$th generation is the
sum of $z_{i}+1$ i.i.d. random variables. Then we have the decomposition
\begin{align}
e^{S_{n}}\mathcal{P}_{(z_{i}+1)\mathbf{e}_{i}}\left( |\mathbf{Z}%
_{n}|=k\right) & =\mathcal{P}_{\mathbf{e}_{i}}\left( |\mathbf{Z}%
_{n}|=0\right) e^{S_{n}}\mathcal{P}_{z_{i}\mathbf{e}_{i}}\left( |\mathbf{Z}%
_{n}|=k\right)  \notag \\
& \quad +\sum_{j=1}^{k-1}\mathcal{P}_{\mathbf{e}_{i}}\left( |\mathbf{Z}%
_{n}|=j\right) e^{S_{n}}\mathcal{P}_{z_{i}\mathbf{e}_{i}}\left( |\mathbf{Z}%
_{n}|=k-j\right)  \notag \\
& \quad +e^{S_{n}}\mathcal{P}_{\mathbf{e}_{i}}\left( |\mathbf{Z}%
_{n}|=k\right) \mathcal{P}_{z_{i}\mathbf{e}_{i}}\left( |\mathbf{Z}%
_{n}|=0\right) \qquad a.s.  \notag
\end{align}%
For the first term, by the assumption (\ref{zi}) and $\mathbb{P}_{\mathbf{e}%
_{i}}\left( |\mathbf{Z}_{n}|=0\right) \overset{n\rightarrow \infty }{%
\longrightarrow }\mathbb{P}_{\mathbf{e}_{i}}\left( |\mathbf{Z}_{\infty
}|=0\right) =\mathcal{R}(i)$ under $\mathbf{P}^{+}$, so
\begin{equation*}
\mathcal{P}_{\mathbf{e}_{i}}\left( |\mathbf{Z}_{n}|=0\right) e^{S_{n}}%
\mathcal{P}_{z_{i}\mathbf{e}_{i}}\left( |\mathbf{Z}_{n}|=k\right) \overset{%
n\rightarrow \infty }{\longrightarrow }\dfrac{z_{i}}{|\mathbf{v}|u_{i}}%
\mathcal{Q}^{2}(i)\mathcal{R}^{z_{i}}(i)\quad \mathbf{P}^{+}-a.s.
\end{equation*}%
For the second term, from previous analysis and the assumption (\ref{zi}),
we know under $\mathbf{P}^{+}$, both $e^{S_{n}}\mathcal{P}_{\mathbf{e}%
_{i}}\left( |\mathbf{Z}_{n}|=j\right) $ and $e^{S_{n}}\mathcal{P}_{z_{i}%
\mathbf{e}_{i},v}\left( |\mathbf{Z}_{n}|=k-j\right) $ converge. Recall that
under $\mathbf{P}^{+}$, $S_{n}\rightarrow \infty $ a.s. So $\mathcal{P}_{%
\mathbf{e}_{i}}\left( |\mathbf{Z}_{n}|=j\right) \overset{n\rightarrow \infty
}{\longrightarrow }0$, $\mathbf{P}^{+}-a.s.$ Thus we get
\begin{equation*}
\sum_{j=1}^{k-1}\mathcal{P}_{\mathbf{e}_{i}}\left( |\mathbf{Z}_{n}|=j\right)
e^{S_{n}}\mathcal{P}_{z_{i}\mathbf{e}_{i}}\left( |\mathbf{Z}_{n}|=k-j\right)
\overset{n\rightarrow \infty }{\longrightarrow }0\quad \mathbf{P}^{+}-a.s.
\end{equation*}%
For the last term, by the independence of the evolution of each particle,
\begin{equation*}
\mathcal{P}_{z_{i}\mathbf{e}_{i}}\left( |\mathbf{Z}_{n}|=0\right) =\mathcal{P%
}_{\mathbf{e}_{i}}\left( |\mathbf{Z}_{n}|=0\right) ^{z_{i}}\overset{%
n\rightarrow \infty }{\longrightarrow }\mathcal{R}^{z_{i}}(i).
\end{equation*}%
Combined with the case of $z_{i}=1$ gives
\begin{equation*}
e^{S_{n}}\mathcal{P}_{\mathbf{e}_{i}}\left( |\mathbf{Z}_{n}|=k\right)
\mathcal{P}_{z_{i}\mathbf{e}_{i}}\left( |\mathbf{Z}_{n}|=0\right) \overset{%
n\rightarrow \infty }{\longrightarrow }\dfrac{1}{|\mathbf{v}|u_{i}}\mathcal{Q%
}^{2}(i)\mathcal{R}^{z_{i}}(i)\quad \mathbf{P}^{+}-a.s.
\end{equation*}%
As a result we obtain
\begin{equation*}
e^{S_{n}}\mathcal{P}_{(z_{i}+1)\mathbf{e}_{i}}\left( |\mathbf{Z}%
_{n}|=k\right) \overset{n\rightarrow \infty }{\longrightarrow }\dfrac{z_{i}+1%
}{|\mathbf{v}|u_{i}}\mathcal{Q}^{2}(i)\mathcal{R}^{z_{i}}(i)\quad \mathbf{P}%
^{+}-a.s.
\end{equation*}%
The first step is accomplished.

\noindent \textit{step 2, induction from $L$-type to $L+1$-type. }

We assume that under $\mathbf{P}^{+}$, for $L\geq 1$ and all $\left(
z_{1},\ldots ,z_{L}\right) \in \mathbb{N}_{0}^{L}$,
\begin{equation*}
e^{S_{n}}\mathcal{P}_{(z_{1},\ldots ,z_{L},0)}\left( |\mathbf{Z}%
_{n}|=k\right) \overset{n\rightarrow \infty }{\longrightarrow }\sum_{i=1}^{L}%
\dfrac{z_{i}}{|\mathbf{v}|u_{i}}\frac{\mathcal{Q}^{2}(i)}{\mathcal{R}(i)}%
\prod\limits_{j=1}^{L}\mathcal{R}^{z_{j}}(j)\quad \mathbf{P}^{+}-a.s.
\end{equation*}%
Then for $\left( z_{1},\ldots ,z_{L+1}\right) \in \mathbb{N}_{0}^{L}$, by
the independence of the evolution of each particle,
\begin{align}
e^{S_{n}}\mathcal{P}_{(z_{1},\ldots ,z_{L+1})}\left( |\mathbf{Z}%
_{n}|=k\right) & =\mathcal{P}_{(z_{1},\ldots ,z_{L},0)}\left( |\mathbf{Z}%
_{n}|=0\right) e^{S_{n}}\mathcal{P}_{z_{L+1}\mathbf{e}_{L+1}}\left( |\mathbf{%
Z}_{n}|=k\right)  \notag \\
& \quad +\sum_{j=1}^{k-1}\mathcal{P}_{(z_{1},\ldots ,z_{L},0)}\left( |%
\mathbf{Z}_{n}|=j\right) e^{S_{n}}\mathcal{P}_{z_{L+1}\mathbf{e}%
_{L+1}}\left( |\mathbf{Z}_{n}|=k-j\right)  \notag \\
& \quad +e^{S_{n}}\mathcal{P}_{(z_{1},\ldots ,z_{L},0)}\left( |\mathbf{Z}%
_{n}|=k\right) \mathcal{P}_{z_{L+1}\mathbf{e}_{L+1}}\left( |\mathbf{Z}%
_{n}|=0\right) .  \notag
\end{align}%
According to step 1 and the assumption, under $\mathbf{P}^{+}$, both $%
e^{S_{n}}\mathcal{P}_{(z_{1},\cdots ,z_{L},0)}\left( |\mathbf{Z}%
_{n}|=j\right) $ and \newline
$e^{S_{n}}\mathcal{P}_{z_{L+1}\mathbf{e}_{L+1}}\left( |\mathbf{Z}%
_{n}|=k-j\right) $ converge to finite limits as $n\rightarrow \infty $. For
the same reason as in step 1, the second term converges to $0$, $\mathbf{P}%
^{+}$-a.s. Similar to the method used in step 1, we get that under $\mathbf{P%
}^{+}$,
\begin{align*}
e^{S_{n}}\mathcal{P}_{(z_{1},\ldots ,z_{L+1})}\left( |\mathbf{Z}%
_{n}|=k\right) & \overset{n\rightarrow \infty }{\longrightarrow }\dfrac{%
z_{L+1}}{|\mathbf{v}|u_{L+1}}\frac{\mathcal{Q}^{2}(L+1)}{\mathcal{R}(L+1)}%
\prod_{j=1}^{L+1}\mathcal{R}^{z_{j}}(j) \\
& +\sum_{i=1}^{L}\dfrac{z_{i}}{|\mathbf{v}|u_{i}}\frac{\mathcal{Q}^{2}(i)}{%
\mathcal{R}(i)}\prod_{j=1}^{L+1}\mathcal{R}^{z_{j}}(j) \\
& \quad =\sum_{i=1}^{L+1}\dfrac{z_{i}}{|\mathbf{v}|u_{i}}\frac{\mathcal{Q}%
^{2}(i)}{\mathcal{R}(i)}\prod_{j=1}^{L+1}\mathcal{R}^{z_{j}}(j).
\end{align*}%
Combining the two steps we conclude that under $\mathbf{P}^{+}$,
\begin{equation}
e^{S_{n}}\mathcal{P}_{(z_{1},\ldots ,z_{K})}\left( |\mathbf{Z}_{n}|=k\right)
\overset{n\rightarrow \infty }{\longrightarrow }\sum_{i=1}^{K}\dfrac{z_{i}}{|%
\mathbf{v}|u_{i}}\frac{\mathcal{Q}^{2}(i)}{\mathcal{R}(i)}\prod_{j=1}^{K}%
\mathcal{R}^{z_{j}}(j).  \label{K}
\end{equation}%
Observe that the limit does not depend on $k$.

By (\ref{K}) and Lemma \ref{dya1}, we complete the proof.\hfill \rule%
{2mm}{3mm}\vspace{4mm}

\begin{corollary}
\label{Cor1}For every $l\in \left\{ 1,\ldots,K\right\} $ and $\mathbf{z}%
=(z_{1},\cdots ,z_{K})\in \mathbb{N}_{0}^{K}$,
\begin{equation*}
T_{l}(\mathbf{z}):=\lim\limits_{n\rightarrow \infty }\mathbf{E}\left[
e^{S_{n}}\mathcal{P}_{\mathbf{z}}\left( \mathbf{Z}_{n}=\mathbf{e}_{l}\right) %
\Big|L_{n}\geq 0\right] =\frac{v_{l}}{\left\vert \mathbf{v}\right\vert }T(%
\mathbf{z}).
\end{equation*}
\end{corollary}

\textbf{Proof. }To prove the corollary it is necessary to take into account (%
\ref{Z_nt}) with $n$ for $n-\lfloor nt\rfloor $, to observe, using formulas (%
\ref{D_indiv}) and (\ref{Dn}), that
\begin{equation}
\frac{D_{n}(l)}{|\mathbf{D}_{n}|}\rightarrow \frac{v_{l}}{|\mathbf{v}|}\quad
\mathbf{P}^{+}-a.s.  \label{Dratiolim}
\end{equation}%
as $n\rightarrow \infty ,$ and to repeat almost literally the proof of Lemma %
\ref{Zj}.

We now formulate the last auxiliary lemma needed to prove Theorem \ref{th2.1}%
.

\begin{lemma}
\textrm{{\textbf{\cite{boinghoff}}}}\label{bo} For every $\varepsilon >0$
and $t\in (0,1)$, there exists an $m\in \mathbb{N}$ such that
\begin{equation*}
\mathbf{P}\left( \tau _{\lfloor nt\rfloor ,n}>\lfloor nt\rfloor +m\Big|%
L_{n}\geq 0\right) \leq \varepsilon .
\end{equation*}
\end{lemma}

\subsection{Proof of Theorem \protect\ref{th2.1}\label{Sub_Proof21}}

First we prove Theorem \ref{th2.1} for the case when $\mathbf{Z}_{n}=\mathbf{%
e}_{l}$ for a fixed $l$.

In the intermediately supercirical case $\mathbf{E}X=0,$ and, therefore,
relation (\ref{Difference2}) holds true implying in view of (\ref{SuperBasic}%
) and (\ref{Repres3}) that
\begin{equation*}
\mathcal{P}_{\mathbf{e}_{i}}\left( \mathbf{Z}_{n}=\mathbf{e}_{l}\right) =%
\frac{D_{n}(l)}{1+|\mathbf{D}_{n}|}\left( \frac{M_{1,n}(i,l)}{D_{n}(l)}-%
\frac{|\mathbf{M}_{1,n}(i)|}{1+|\mathbf{D}_{n}|}\right) =\frac{D_{n}(l)}{|%
\mathbf{D}_{n}|}\frac{|\mathbf{M}_{1,n}(i)|}{(1+|\mathbf{D}_{n}|)^{2}}%
(1+\varepsilon _{i}(n)),
\end{equation*}%
where $\varepsilon _{i}(n)\rightarrow 0$ $\mathbf{P}$-a.s as $n\rightarrow
\infty $. After the change of measure we see that%
\begin{equation}
\mathbb{P}_{\mathbf{e}_{i}}(\mathbf{Z}_{n}=\mathbf{e}_{j})\sim \kappa ^{n}%
\mathbf{E}\left[ \frac{e^{S_{n}}}{\left\vert \mathbf{M}_{1,n}(i)\right\vert }%
\frac{D_{n}(l)}{|\mathbf{D}_{n}|}\mathcal{Q}_{n}^{2}\left( i\right) \right] ,
\label{10}
\end{equation}%
\ as $n\rightarrow \infty $. Put
\begin{equation*}
V_{n}:=\frac{e^{S_{n}}}{\left\vert \mathbf{M}_{1,n}(i)\right\vert }\frac{%
D_{n}(l)}{|\mathbf{D}_{n}|}\mathcal{Q}_{n}\left( i\right) \cdot \mathbbm{1}%
_{\{\mathbf{Z}_{0}=\mathbf{e}_{i},|\mathbf{Z}_{n}|>0\}}.
\end{equation*}

For $k\leq n$ the event $\{|\mathbf{Z}_{n}|>0\}$ implies $\{|\mathbf{Z}%
_{k}|>0\}$. Thus,
\begin{align}
V_{k,n}:& =\mathbf{E}\left[ V_{n};|\mathbf{Z}_{k}|>0\Big|L_{k,n}\geq 0,%
\mathcal{F}_{k}\right]  \notag \\
& =\mathbf{E}\left[ \frac{e^{S_{n}}}{\left\vert \mathbf{M}%
_{1,n}(i)\right\vert }\frac{D_{n}(l)}{|\mathbf{D}_{n}|}\cdot \mathbbm{1}_{\{%
\mathbf{Z}_{0}=\mathbf{e}_{i},|\mathbf{Z}_{n}|>0\}};|\mathbf{Z}_{k}|>0\Big|%
L_{k,n}\geq 0,\mathcal{F}_{k}\right]  \notag \\
& =\mathbf{E}\left[ \frac{e^{S_{n}}}{\left\vert \mathbf{M}%
_{1,n}(i)\right\vert }\frac{D_{n}(l)}{|\mathbf{D}_{n}|}\cdot \mathbbm{1}_{\{%
\mathbf{Z}_{0}=\mathbf{e}_{i},|\mathbf{Z}_{n}|>0\}}\Big|L_{k,n}\geq 0,%
\mathcal{F}_{k}\right] \quad a.s.  \notag
\end{align}%
Note that $V_{k,n}$ is adapted to the flow of $\sigma $-algebras $\{\mathcal{%
F}_{k}\}_{k\in \mathbb{N}}$. Denoting $\mathbf{P}_{\mu _{k}}^{+}$ the
measure corresponding to $\mathbf{P}^{+}$ and calculated given $\mathcal{F}%
_{k}$ , we see by (\ref{Mi}) and Lemma \ref{u} that for some random vector $%
\mathbf{u}=\left( u_{1},\ldots ,u_{K}\right) $ and all $i=1,2,\ldots ,K$
\begin{equation}
e^{-S_{n}}\left\vert \mathbf{M}_{1,n}(i)\right\vert \sim |\mathbf{v}|u_{i}(%
\mathbf{M}_{1,n})=|\mathbf{v}|u_{i}(\mathbf{M}_{1,k}\mathbf{M}_{k+1,n})%
\overset{n\rightarrow \infty }{\longrightarrow }|\mathbf{v}|u_{i}\qquad
\mathbf{P}_{\mu _{k}}^{+}-a.s.  \label{Ms}
\end{equation}%
and, in addition,%
\begin{equation*}
\mathbbm{1}_{\{\mathbf{Z}_{0}=\mathbf{e}_{i},|\mathbf{Z}_{n}|>0\}}\overset{%
n\rightarrow \infty }{\longrightarrow }\mathbbm{1}_{\{\mathbf{Z}_{0}=\mathbf{%
e}_{i},|\mathbf{Z}_{\infty }|>0\}}\qquad \mathbf{P}_{\mu _{k}}^{+}-a.s.,
\end{equation*}%
\begin{equation}
\mathcal{Q}_{n}\left( i\right) =1-F_{0,k}^{(i)}(\mathbf{F}_{k+1,n}(\mathbf{0}%
))\overset{n\rightarrow \infty }{\longrightarrow }1-F_{0,k}^{(i)}(\mathbf{F}%
_{k+1,\infty }(\mathbf{0}))\qquad \mathbf{P}_{\mu _{k}}^{+}-a.s.
\label{DefQ-limit}
\end{equation}%
Recalling (\ref{Dratiolim}) we see that the conditions of Lemma \ref{dya1}
are valid for $V_{k,n}$. As a result we get that $\mathbf{P}_{\mu
_{k}}^{+}-a.s.$, as $n\rightarrow \infty $
\begin{equation}
V_{k,n}\rightarrow \mathbf{E}_{\mu _{k}}^{+}\left[ \frac{1-F_{0,k}^{(i)}(%
\mathbf{F}_{k+1,\infty }(\mathbf{0}))}{u_{i}}\mathbbm{1}_{\{\mathbf{Z}_{0}=%
\mathbf{e}_{i},|\mathbf{Z}_{\infty }|>0\}}\right] :=V_{k,\infty }^{(i)}.
\label{DefVkn}
\end{equation}%
Since $\mathbf{P}(L_{n}\geq 0)\sim \mathbf{P}(L_{k,n}\geq 0)$ as $%
n\rightarrow \infty $ for each fixed $k$, we get
\begin{align}
\mathbf{E}\left[ V_{n};|\mathbf{Z}_{k}|>0,L_{k,n}\geq 0\Big|\mathcal{F}_{k}%
\right] & =\mathbf{P}(L_{k,n}\geq 0)\mathbf{E}\left[ V_{k,n}\Big|L_{k,n}\geq
0,\mathcal{F}_{k}\right]  \notag \\
& =\mathbf{P}(L_{n}\geq 0)\left( V_{k,\infty }^{(i)}+o(1)\right) \ \quad
\mathbf{P}_{\mu _{k}}^{+}-a.s.  \notag
\end{align}%
This representaion allows us to use Lemma \ref{dya2} to obtain
\begin{equation}
\mathbf{E}\left[ V_{n};|\mathbf{Z}_{\tau _{n}}|>0\right] =\frac{v_{l}}{|%
\mathbf{v}|^{2}}\left(\hat{\Delta} _{i}(\mathbf{e}_{l})+o(1)\right) \mathbf{P}%
(L_{n}\geq 0)
\end{equation}
where
\begin{equation}
\hat{\Delta}_{i}(\mathbf{e}_{l}):=\sum\limits_{k=0}^{\infty }\mathbf{E}\left[
V_{k,\infty }^{(i)};\tau _{k}=k\right] .  \label{DefDelta_i}
\end{equation}%
According to Lemma \ref{dya2} $\hat{\Delta} _{i}(\mathbf{e}_{l})<\infty .$
Besides, $\hat{\Delta} _{i}(\mathbf{e}_{l})\geq \mathbf{E}\left[ V_{0,\infty
}^{(i)};\tau _{0}=0\right] >0$. Observing that
\begin{align}
\mathbf{E}\left[ V_{n};|\mathbf{Z}_{\tau _{n}}|>0\right] & =\mathbf{E}\left[
\frac{e^{S_{n}}}{\left\vert \mathbf{M}_{1,n}(i)\right\vert }\frac{D_{n}(l)}{|%
\mathbf{D}_{n}|}\mathcal{Q}_{n}\left( i\right) \cdot \mathbbm{1}_{\{\mathbf{Z%
}_{0}=\mathbf{e}_{i},|\mathbf{Z}_{n}|>0\}};|\mathbf{Z}_{\tau _{n}}|>0\right]
\notag \\
& =\mathbf{E}\left[ \frac{e^{S_{n}}}{\left\vert \mathbf{M}%
_{1,n}(i)\right\vert }\frac{D_{n}(l)}{|\mathbf{D}_{n}|}\mathcal{Q}_{n}\left(
i\right) \cdot \mathbbm{1}_{\{\mathbf{Z}_{0}=\mathbf{e}_{i},|\mathbf{Z}%
_{n}|>0\}}\right]  \notag \\
& =\mathbf{E}\left[ \frac{e^{S_{n}}}{\left\vert \mathbf{M}%
_{1,n}(i)\right\vert }\frac{D_{n}(l)}{|\mathbf{D}_{n}|}\mathcal{Q}%
_{n}^{2}\left( i\right) \right]  \label{V-gener}
\end{align}%
and applying (\ref{10}) and Lemma 2.1 in \cite{AGKV2005}, we see for $\Delta
_{i}(\mathbf{e}_{l}):=|\mathbf{v}|^{2}\hat{\Delta}_{i}(\mathbf{e}_{l})/v_{l}$
that, as $n\rightarrow \infty $,
\begin{equation*}
\mathbb{P}_{\mathbf{e}_{i}}(\mathbf{Z}_{n}=\mathbf{e}_{l})\sim \kappa ^{n}%
\Delta _{i}(\mathbf{e}_{l})\mathbf{P}%
(L_{n}\geq 0)\sim \kappa ^{n}\Delta _{i}(%
\mathbf{e}_{l})\cdot n^{-(1-a)}l(n),
\end{equation*}
where $l(n)$ is a function slowly varying at infinity. \ Note, finally, that
in view of (\ref{V-gener}) $\Delta _{i}(\mathbf{e}_{l})$ may be rewritten as%
\begin{equation}
\Delta _{i}(\mathbf{e}_{l})=\sum\limits_{k=0}^{\infty }\mathbf{E}\left[
\frac{\left( 1-F_{0,k}^{(i)}(\mathbf{F}_{k+1,\infty }(\mathbf{0}))\right)
^{2}}{u_{i}};\tau _{k}=k\right]  \label{Final_delta}
\end{equation}

Using now methods similar to those applied to prove Theorem \ref{th1.1}, one
can check the validity of Theorem \ref{th2.1} for any nonzero vector $%
\mathbf{z}\in \mathbb{N}^{K}$ with%
\begin{equation*}
\Delta _{i}(\mathbf{z})=C(\mathbf{z},\mathbf{v})\Delta _{i}(\mathbf{e}_{l}).
\end{equation*}

\hfill\rule{2mm}{3mm}\vspace{4mm}

\subsection{Proof of Theorem \protect\ref{T_uniformInterm}}

Recall the definition of the measure $\mathbf{P}$ and decompose at the
minimum, then for every $0\leq m\leq n$ and every $k\in \mathbb{N}$,
\begin{align}
\mathbb{P}_{\mathbf{e}_{i}}\left( |\mathbf{Z}_{n}|=k\right) & =\sum_{j=0}^{m}%
\mathbb{P}_{\mathbf{e}_{i}}\left( |\mathbf{Z}_{n}|=k;\tau _{n}=j\right) +%
\mathbb{P}_{\mathbf{e}_{i}}\left( |\mathbf{Z}_{n}|=k;\tau _{n}>m\right)
\notag \\
& =\sum_{j=0}^{m}\mathbb{P}_{\mathbf{e}_{i}}\left( |\mathbf{Z}_{n}|=k;\tau
_{j}=j,L_{j,n}\geq 0\right) +\mathbb{P}_{\mathbf{e}_{i}}\left( |\mathbf{Z}%
_{n}|=k;\tau _{n}>m\right)  \notag \\
& =\kappa ^{n}\sum_{j=0}^{m}\mathbf{E}\left[ e^{S_{n}}\mathcal{P}_{\mathbf{e}%
_{i}}\left( |\mathbf{Z}_{n}|=k\right) ;\tau _{j}=j,L_{j,n}\geq 0\right] +%
\mathbb{P}_{\mathbf{e}_{i}}\left( |\mathbf{Z}_{n}|=k;\tau _{n}>m\right) .
\notag
\end{align}%
For the first term, by the Markov property of the process, Lemma \ref{Zj}
and the dominated convergence theorem we have
\begin{align}
& \quad \ \kappa ^{n}\sum_{j=0}^{m}\mathbf{E}\left[ e^{S_{n}}\mathcal{P}_{%
\mathbf{e}_{i}}\left( |\mathbf{Z}_{n}|=k\right) ;\tau _{j}=j,L_{j,n}\geq 0%
\right]  \notag \\
& =\kappa ^{n}\sum_{j=0}^{m}\mathbf{E}\left[ e^{S_{j}}\mathbf{E}\left[
e^{S_{n-j}}\mathcal{P}_{\mathbf{Z}_{j}}\left( |\mathbf{Z}_{n}|=k\right) \Big|%
L_{n-j}\geq 0\right] ;\tau _{j}=j\right] \mathbf{P}(L_{n-j}\geq 0)  \notag \\
& \sim \kappa ^{n}\sum_{j=0}^{m}\mathbf{E}\left[ e^{S_{j}}T(\mathbf{Z}%
_{j});\tau _{j}=j\right] \mathbf{P}(L_{n-j}\geq 0)
\end{align}%
as $n\rightarrow \infty .$ Hence this sum does not depend on $k$ as $%
n\rightarrow \infty $.

For the second term, by Lemma \ref{l7} and Theorem \ref{th2.1}, we know that
for every $\varepsilon >0$ and $m$ big enough, as $n\rightarrow \infty $,
\begin{equation*}
\mathbb{P}_{\mathbf{e}_{i}}\left( |\mathbf{Z}_{n}|=k;\tau _{n}>m\right) \leq
\varepsilon \mathbb{P}_{\mathbf{e}_{i}}\left( |\mathbf{Z}_{n}|=1\right) .
\end{equation*}%
Hence, for every $k\in \mathbb{N}$,
\begin{equation}
\dfrac{\lim\limits_{n\rightarrow \infty }\mathbb{P}_{\mathbf{e}_{i}}\left( |%
\mathbf{Z}_{n}|=1\right) }{\lim\limits_{n\rightarrow \infty }\mathbb{P}_{%
\mathbf{e}_{i}}\left( |\mathbf{Z}_{n}|=k\right) }=1  \notag
\end{equation}%
which, evidently, imply the statement of Theorem \ref{T_uniformInterm}. \ \
\ \ \ \ \ \ \ \ \ \ \ \ \ \ \ \ \ \ \ \ \ \ \ \ \ \ \ \ \ \ \ \ \ \ \ \ \ \
\ \ \ \ \ \rule{2mm}{3mm}\vspace{4mm} \ \ \ \ \ \ \ \ \ \ \ \ \ \ \ \ \ \ \
\ \ \ \ \ \ \ \ \ \ \ \ \ \ \ \ \ \ \ \ \ \ \ \ \ \ \ \ \ \ \ \ \ \ \ \
\hfill

\subsection{Proof of Theorem \protect\ref{T_cond_interm}}

For $l\in \left\{ 1,\ldots ,K\right\} ,0\leq m\leq n$ and $t\in (0,1)$ write
\begin{equation*}
\mathbb{P}_{\mathbf{e}_{i}}\left( \mathbf{Z}_{\tau _{\lfloor nt\rfloor ,n}}=%
\mathbf{z},\mathbf{Z}_{n}=\mathbf{e}_{l}\right) =\mathbb{P}_{\mathbf{e}%
_{i}}\left( \mathbf{Z}_{\tau _{\lfloor nt\rfloor ,n}}=\mathbf{z},\mathbf{Z}%
_{n}=\mathbf{e}_{l};\tau _{n}\leq m\right) +\mathbb{P}_{\mathbf{e}%
_{i}}\left( \mathbf{Z}_{\tau _{\lfloor nt\rfloor ,n}}=\mathbf{z},\mathbf{Z}%
_{n}=\mathbf{e}_{l};\tau _{n}>m\right) .
\end{equation*}
By Lemma \ref{l7} and Theorem \ref{th2.1} the second term admits the
following estimate for every $\varepsilon >0$ and $m$ big enough:
\begin{equation*}
\limsup_{n\rightarrow \infty }\dfrac{\mathbb{P}_{\mathbf{e}_{i}}\left(
\mathbf{Z}_{\tau _{\lfloor nt\rfloor ,n}}=\mathbf{z},\mathbf{Z}_{n}=\mathbf{e%
}_{l};\tau _{n}>m\right) }{\mathbb{P}_{\mathbf{e}_{i}}\left( \mathbf{Z}_{n}=%
\mathbf{e}_{l}\right) }\leq \varepsilon .
\end{equation*}%
So this term can be neglected when we calculate the conditional probability,
and we only need to consider the first term.

For fixed $k\leq m$ and $r<n-\lfloor nt\rfloor $, using the change of
measure, we have
\begin{align}
& \quad \mathbb{P}_{\mathbf{e}_{i}}\left( \mathbf{Z}_{\tau _{\lfloor
nt\rfloor ,n}}=\mathbf{z},\mathbf{Z}_{n}=\mathbf{e}_{l};\tau _{n}=k\right) =%
\mathbb{P}_{\mathbf{e}_{i}}\left( \mathbf{Z}_{\tau _{\lfloor nt\rfloor ,n}}=%
\mathbf{z},\mathbf{Z}_{n}=\mathbf{e}_{l};\tau _{k}=k,L_{k,n}\geq 0\right)
\notag \\
& =\kappa ^{n}\sum_{j=0}^{n-\lfloor nt\rfloor }\mathbf{E}\left[ e^{S_{n}}%
\mathcal{P}_{\mathbf{e}_{i}}\left( \mathbf{Z}_{\tau _{\lfloor nt\rfloor ,n}}=%
\mathbf{z},\mathbf{Z}_{n}=\mathbf{e}_{l}\right) ;\tau _{k}=k,L_{k,n}\geq
0,\tau _{\lfloor nt\rfloor ,n}=\lfloor nt\rfloor +j\right]   \notag \\
& =\kappa ^{n}\left( H_{1}(k,r)+H_{2}(k,r)\right) ,  \label{slt}
\end{align}%
where
\begin{equation*}
H_{1}(k,r):=\sum\limits_{j=0}^{r}\mathbf{E}\left[ e^{S_{\lfloor nt\rfloor
+j}}e^{S_{n}-S_{\lfloor nt\rfloor +j}}\mathcal{P}_{\mathbf{e}_{i}}\left(
\mathbf{Z}_{_{\lfloor nt\rfloor +j}}=\mathbf{z},\mathbf{Z}_{n}=\mathbf{e}%
_{l}\right) ;\tau _{k}=k,L_{k,n}\geq 0,\tau _{\lfloor nt\rfloor ,n}=\lfloor
nt\rfloor +j\right]
\end{equation*}%
and \
\begin{equation}
H_{2}(k,r):=\mathbf{E}\left[ e^{S_{n}}\mathcal{P}_{\mathbf{e}_{i}}\left(
\mathbf{Z}_{_{\lfloor nt\rfloor +j}}=\mathbf{z},\mathbf{Z}_{n}=\mathbf{e}%
_{l}\right) ;\tau _{k}=k,L_{k,n}\geq 0,\tau _{\lfloor nt\rfloor ,n}>\lfloor
nt\rfloor +r\right] .  \label{DefH2}
\end{equation}%
Observe that $\{\tau _{\lfloor nt\rfloor ,n}=j\}=\{\tau _{\lfloor nt\rfloor
,j}=j\}\cap \{L_{j,n}\geq 0\}$ and
\begin{equation*}
\{\tau _{\lfloor nt\rfloor ,n}=j,L_{n}\geq 0\}=\{\tau _{\lfloor nt\rfloor
,j}=j,L_{j}\geq 0\}\cap \{L_{j,n}\geq 0\}.
\end{equation*}%
Hence
\begin{align}
H_{1}(k,r)=& \sum\limits_{j=0}^{r}\mathbf{E}\left[ e^{S_{\lfloor nt\rfloor
+j}}\mathcal{P}_{\mathbf{e}_{i}}\left( \mathbf{Z}_{\lfloor nt\rfloor +j}=%
\mathbf{z}\right) ;\tau _{k}=k,L_{k,\lfloor nt\rfloor +j}\geq 0,\tau
_{\lfloor nt\rfloor ,\lfloor nt\rfloor +j}=\lfloor nt\rfloor +j\right]
\notag \\
& \times \mathbf{E}\left[ e^{S_{n-\lfloor nt\rfloor -j}}\mathcal{P}_{\mathbf{%
z}}\left( \mathbf{Z}_{n-\lfloor nt\rfloor -j}=\mathbf{e}_{l}\right)
;L_{n-\lfloor nt\rfloor -j}\geq 0\right]   \notag \\
=& \sum\limits_{j=0}^{r}\mathbf{E}\left[ e^{S_{\lfloor nt\rfloor +j}}%
\mathcal{P}_{\mathbf{e}_{i}}\left( \mathbf{Z}_{\lfloor nt\rfloor +j}=\mathbf{%
z}\right) ;\tau _{k}=k,L_{k,\lfloor nt\rfloor +j}\geq 0,\tau _{\lfloor
nt\rfloor ,\lfloor nt\rfloor +j}=\lfloor nt\rfloor +j\right]   \notag \\
& \times \mathbf{E}\left[ e^{S_{n-\lfloor nt\rfloor -j}}\mathcal{P}_{\mathbf{%
z}}\left( \mathbf{Z}_{n-\lfloor nt\rfloor -j}=\mathbf{e}_{l}\right) \Big|%
L_{n-\lfloor nt\rfloor -j}\geq 0\right] \mathbf{P}(L_{n-\lfloor nt\rfloor
-j}\geq 0).  \label{s1a}
\end{align}%
Recalling that $\mathbf{E}X=0$ and using (\ref{SuperBasic}), (\ref%
{Difference2}), (\ref{Ms}) and (\ref{Dratiolim}) we conclude that, as $%
n\rightarrow \infty $
\begin{equation*}
e^{S_{\lfloor nt\rfloor +j}}\mathcal{P}_{\mathbf{e}_{i}}\left( \mathbf{Z}%
_{\lfloor nt\rfloor +j}=\mathbf{z}\right) -W_{nt,j}(\mathbf{z})=o(1)\quad
\mathbf{P-}a.s.,
\end{equation*}%
where%
\begin{equation}
W_{nt,j}(\mathbf{z}):=C(\mathbf{z},\mathbf{v})\dfrac{\mathcal{Q}_{\lfloor
nt\rfloor +j}^{2}(i)}{u_{i}(M_{1,\lfloor nt\rfloor +j})}.  \label{DefW}
\end{equation}%
By Corollary \ref{Cor1} we have
\begin{equation}
\mathbf{E}\left[ e^{S_{n-\lfloor nt\rfloor -j}}\mathcal{P}_{\mathbf{z}%
}\left( \mathbf{Z}_{n-\lfloor nt\rfloor -j}=\mathbf{e}_{l}\right) \Big|%
L_{n-\lfloor nt\rfloor -j}\geq 0\right] \rightarrow T_{l}(\mathbf{z})\quad
\mathbf{P-}a.s.  \notag
\end{equation}%
as $n-\lfloor nt\rfloor -j\rightarrow \infty $. Inserting the above formulas
into (\ref{s1a}) gives
\begin{eqnarray*}
H_{1}(k,r) &\sim &\sum\limits_{j=0}^{r}\mathbf{E}\Big[W_{nt,j}(\mathbf{z}%
);\tau _{k}=k,L_{k,\lfloor nt\rfloor +j}\geq 0,\tau _{\lfloor nt\rfloor
,\lfloor nt\rfloor +j}=\lfloor nt\rfloor +j\Big] \\
&&\qquad \times T_{l}(\mathbf{z})\mathbf{P}(L_{n-\lfloor nt\rfloor -j}\geq 0)
\\
&=&\sum\limits_{j=0}^{r}\mathbf{E}\Big[W_{nt,j}(\mathbf{z});\tau
_{k}=k,L_{k,n}\geq 0,\tau _{\lfloor nt\rfloor ,n}=\lfloor nt\rfloor +j\Big]%
T_{l}(\mathbf{z}) \\
&=&\mathbf{E}\Big[W_{nt,j}(\mathbf{z});\tau _{k}=k,L_{k,n}\geq 0,\tau
_{\lfloor nt\rfloor ,n}\leq \lfloor nt\rfloor +r\Big]T_{l}(\mathbf{z}).
\end{eqnarray*}%
Note that by Lemma \ref{bo} and the asymptotic representation $\mathbf{P}%
\left( L_{n}\geq 0\right) \sim n^{-(1-a)}l(n)$ as $n\rightarrow \infty $
(see, for instance, Lemma 2.1 in \cite{AGKV2005}), we have
\begin{align*}
& \mathbf{E}\left[ W_{nt,j}(\mathbf{z});\tau _{k}=k,L_{k,n}\geq 0,\tau
_{\lfloor nt\rfloor ,n}\leq \lfloor nt\rfloor +r\right]  \\
& \qquad =\mathbf{E}\left[ \mathbf{E}\left[ W_{nt,j}(\mathbf{z});\tau
_{\lfloor nt\rfloor ,n}\leq \lfloor nt\rfloor +r\Big|L_{k,n}\geq 0,\mathcal{F%
}_{k}\right] ;\tau _{k}=k\right] \mathbf{P}(L_{n-k}\geq 0) \\
& \qquad \sim \mathbf{E}\left[ \mathbf{E}\left[ W_{nt,j}(\mathbf{z})\Big|%
L_{k,n}\geq 0,\mathcal{F}_{k}\right] ;\tau _{k}=k\right] \mathbf{P}%
(L_{n}\geq 0) \\
& \qquad =C(\mathbf{z},\mathbf{v})\mathbf{E}\left[ \mathbf{E}\left[ \dfrac{%
\mathcal{Q}_{\lfloor nt\rfloor +j}^{2}(i)}{u_{i}(M_{1,\lfloor nt\rfloor +j})}%
\Big|L_{k,n}\geq 0,\mathcal{F}_{k}\right] ;\tau _{k}=k\right] \mathbf{P}%
(L_{n}\geq 0)
\end{align*}%
as $n\rightarrow \infty $. Using (\ref{DefQ-limit}) $\ $we deduce by the
dominated convergence theorem and (\ref{Ms}) that $\mathbf{P}_{\mu
_{k}}^{+}-a.s.$, as $n\rightarrow \infty $%
\begin{equation*}
\mathbf{E}\left[ \dfrac{\mathcal{Q}_{\lfloor nt\rfloor +j}^{2}(i)}{%
u_{i}(M_{1,\lfloor nt\rfloor +j})}\Big|L_{k,n}\geq 0,\mathcal{F}_{k}\right]
\sim \mathbf{E}_{\mu _{k}}^{+}\left[ \dfrac{\left( 1-F_{0,k}^{(i)}(\mathbf{F}%
_{k+1,\infty }(\mathbf{0}))\right) ^{2}}{u_{i}}\right] .
\end{equation*}

By Theorem \ref{th2.1}, Lemma \ref{dya1} and Lemma \ref{bo}, we get if first
$n\rightarrow \infty ,$ than $r\rightarrow \infty $
\begin{align}
& \mathbf{E}\left[ \mathbf{E}\left[ W_{nt,j}(\mathbf{z});\tau _{\lfloor
nt\rfloor ,n}\leq \lfloor nt\rfloor +r\Big|L_{k,n}\geq 0,\mathcal{F}_{k}%
\right] ;\tau _{k}=k\right]   \notag \\
& \rightarrow C(\mathbf{z},\mathbf{v})\mathbf{E}\left[ \mathbf{E}_{\mu
_{k}}^{+}\left[ \dfrac{\left( 1-F_{0,k}^{(i)}(\mathbf{F}_{k+1,\infty }(%
\mathbf{0}))\right) ^{2}}{u_{i}}\right] ;\tau _{k}=k\right] :=\psi _{ik}(%
\mathbf{z}).  \notag
\end{align}%
Combining the estimates above we conclude that
\begin{equation}
\mathbb{P}_{\mathbf{e}_{i}}\left( \mathbf{Z}_{\tau _{\lfloor nt\rfloor ,n}}=%
\mathbf{z},\mathbf{Z}_{n}=\mathbf{e}_{l};\tau _{n}=k\right) =\psi _{ik}(%
\mathbf{z})T_{l}(\mathbf{z})\mathbf{P}(L_{n}\geq 0)(1+\varepsilon (n,k)),
\label{SingTerm}
\end{equation}%
where $\lim_{n\rightarrow \infty }\varepsilon (n,k)=0.$

To evaluate $H_{2}(k,r)$ specified in (\ref{DefH2}) we first recall that by (%
\ref{z}) and (\ref{Bound_M})
\begin{equation*}
\mathcal{P}_{\mathbf{e}_{i}}\left( \mathbf{Z}_{\lfloor nt\rfloor +j}=\mathbf{%
z},\mathbf{Z}_{n}=\mathbf{e}_{l}\right) \leq \mathcal{P}_{\mathbf{e}%
_{i}}\left( |\mathbf{Z}_{n}|=1\right) \leq \frac{\left\vert \mathbf{v}%
\right\vert }{v_{i}}e^{-S_{n}}\quad a.s.
\end{equation*}%
Using this estimate, monotonicity of $\mathbf{P}\left( L_{n}\geq 0\right) $
in $n,$ Lemma \ref{bo} and Theorem \ref{th2.1} we get for each fixed $k$ and
sufficiently large $n$
\begin{align}
H_{2}(k,r)& \leq \frac{\left\vert \mathbf{v}\right\vert }{v_{i}}%
\sum_{j=r+1}^{n-\lfloor nt\rfloor }\mathbf{P}\left( \tau _{k}=k,L_{k,n}\geq
0,\tau _{\lfloor nt\rfloor ,n}=\lfloor nt\rfloor +j\right)   \notag \\
& =\frac{\left\vert \mathbf{v}\right\vert }{v_{i}}\mathbf{P}\left( \tau
_{k}=k\right) \mathbf{P}\left( L_{n-k}\geq 0\right) \mathbf{P}\left( \tau
_{\lfloor nt\rfloor -k,n-k}>\lfloor nt\rfloor +r-k\Big|L_{n-k}\geq 0\right)
\notag \\
& \leq \frac{\varepsilon \left\vert \mathbf{v}\right\vert }{v_{i}}\mathbf{P}%
\left( \tau _{k}=k\right) \mathbf{P}\left( L_{n}\geq 0\right) \leq
\varepsilon C\kappa ^{-n}\mathbf{P}\left( \tau _{k}=k\right) \mathbb{P}_{%
\mathbf{e}_{i}}\left( \mathbf{Z}_{n}=\mathbf{e}_{l}\right) .  \label{H2estim}
\end{align}%
Hence when consider the conditional probability, as $n\rightarrow \infty $
and $m\rightarrow \infty $ the sum
\begin{equation*}
\sum\limits_{k=0}^{m}H_{2}(k,r)\leq \varepsilon C\kappa ^{-n}\mathbb{P}_{%
\mathbf{e}_{i}}\left( \mathbf{Z}_{n}=\mathbf{e}_{l}\right)
\end{equation*}%
can be neglected. \ Combining this estimate with (\ref{SingTerm}) gives, as $%
n\rightarrow \infty $
\begin{equation}
\mathbb{P}_{\mathbf{e}_{i}}\left( \mathbf{Z}_{\tau _{\lfloor nt\rfloor ,n}}=%
\mathbf{z},\mathbf{Z}_{n}=\mathbf{e}_{l}\right) \sim \kappa ^{n}C(\mathbf{z},%
\mathbf{v})\Psi _{i}(\mathbf{z})T_{l}(\mathbf{z})\mathbf{P}(L_{n}\geq 0),
\notag
\end{equation}%
where%
\begin{equation*}
\Psi _{i}(\mathbf{z})=\sum_{k=0}^{\infty }\psi _{ik}(\mathbf{z})=C(\mathbf{z}%
,\mathbf{v})\sum_{k=0}^{\infty }\mathbf{E}\left[ \mathbf{E}_{\mu _{k}}^{+}%
\left[ \dfrac{\left( 1-F_{0,k}^{(i)}(\mathbf{F}_{k+1,\infty }(\mathbf{0}%
))\right) ^{2}}{u_{i}}\right] ;\tau _{k}=k\right] .
\end{equation*}%
Using the conclusions of Theorem \ref{th2.1} and (\ref{Final_delta}), we see
that
\begin{equation}
\lim_{n\rightarrow \infty }\mathbb{P}_{\mathbf{e}_{i}}\left( \mathbf{Z}%
_{\tau _{\lfloor nt\rfloor ,n}}=\mathbf{z}\Big|\mathbf{Z}_{n}=\mathbf{e}%
_{l}\right) =C(\mathbf{z},\mathbf{v})T(\mathbf{z})=:q(\mathbf{z}).  \notag
\end{equation}%
By similar method as the proof of Theorem \ref{th1.1}, one can check that $q(%
\mathbf{z})$ is a probability measure on $\mathbb{N}_{0}^{K}$. \hfill \rule%
{2mm}{3mm}\vspace{4mm}

\hfill

\end{document}